\documentclass{amsart}

\usepackage{amsmath}
\usepackage{physics}
\usepackage{eurosym}
\usepackage{amsfonts}
\usepackage{amssymb}
\usepackage{amsmath}
\usepackage{amsfonts}
\usepackage{graphicx}
\usepackage{subfig}
\usepackage{amssymb}
\usepackage{amsmath}
\usepackage{amsfonts}
\usepackage{graphicx}
\usepackage{subfig}
\usepackage{amsfonts}

\theoremstyle{plain}

\newtheorem{definition}{Definition}

\numberwithin{equation}{section}

\begin{document}
\title[]{Notes on the Geometry of Electromagnetic Fields and Maxwell's
Equations along a non-null curves in non flat-3D space forms $M_{q}^{3}(c)$}
\author{Fatma Almaz}
\address{department of mathematics, faculty of arts and sciences, batman
university, batman/ t\"{u}rk\.{ı}ye}
\email{fatma.almaz@batman.edu.tr}
\author{Cumal\.{ı} Ek\.{ı}c\.{ı}}
\address{department of mathematics and computer sciences, faculty of
sciences, esk\.{ı}\c{s}eh\.{ı}r osmangazi university, esk\.{ı}\c{s}eh\.{ı}r/ t\"{u}rk\.{ı}ye}
\email{cekici@ogu.edu.tr}
\thanks{This paper is in final form and no version of it will be submitted
for publication elsewhere.}

\begin{abstract}
In this paper, the directional derivatives in accordance with the
orthonormal frame $\left\{ T,N,B\right\} $ are defined in $M_{q}^{3}(c)$,
the extended Serret-Frenet relations by using Frenet formulas are expressed.
Furthermore, we express the bending elastic energy function for the same
particle in $M_{q}^{3}(c)$ according to curve $\alpha (s,\xi ,\eta )$ and
geometrical interpretation of the energy for unit vector fields and we also
solve Maxwell's equations for the electric and magnetic field vectors in $%
M_{q}^{3}(c).$
\end{abstract}

\maketitle

\section{Introduction}

Maxwell's equations, one of the most elegant and powerful sets of equations
in physics, unified the field of classical electromagnetism by revealing the
deep connection between electricity and magnetism. These four equations,
formulated by James Clerk Maxwell, comprehensively explain how electric
charges and currents generate electric and magnetic fields and how these
fields propagate through space and time. Maxwell's equations not only
explained static electric and magnetic phenomena but also predicted that
changing electric fields could induce magnetic fields, and that changing
magnetic fields could induce electric fields. This interplay theoretically
established the existence of electromagnetic waves propagating at the speed
of light and established the unification between optics and
electromagnetism. These equations form the basis of many technologies we use
in our daily lives, such as radio, television, cell phones, fiber optic
communication, and electric motors. Understanding Maxwell's equations is key
to understanding our electromagnetic universe.

The energy of an electromagnetic field is carried by the field itself and is
expressed in terms of the electric field vector and the magnetic field
vector. Maxwell's equations describe how these fields exist and interact,
while the energy attributed to the fields and the flow of that energy are
also described through these vector fields.

The spaces $M_{q}^{3}(c)$ generalize to 3-dimensional spaces that differ
from Euclidean geometry (curved) and whose metric structure is
pseudo-Riemann. The values of $q$ and $c$ determine the geometric and
topological properties of the space. The physical interpretation of the
space forms $M_{q}^{3}(c)$ depends on the nature of the physical theory
defined within it (e.g., electromagnetism, gravity) and the metric structure
of space. The most common physical interpretation is that, in the case $q=1$%
, these spaces represent possible spacetime solutions in general relativity.
The curvature $(c)$ indicates how this spacetime is bent or curved (the
effect of gravity), while the sign convention $(q)$ determines how the
spacetime and time dimensions are separated and the structure of causality.
Non-flat cases $(c\neq 0)$ are important in scenarios where gravity is
important or in theoretical investigations where the topology and geometry
of spacetime influence physical phenomena. To summarize the physical
significance of the space forms $M_{q}^{3}(c)$ depends strongly on the sign
convention of their metrics $(q)$ and their curvature $(c)$. In the
Lorentzian case $(q=1)$, these spaces are crucial for describing the
geometry of spacetime in gravitational and cosmological models. In the
Riemann case $(q=0)$, they can play a role as configuration spaces or in
theoretical physics models. The constant curvature makes these spaces easier
to treat mathematically, and their high symmetry makes them fertile ground
for research in theoretical physics.

This study examines the fundamental nature of electromagnetic fields,
particularly the interactions between electric and magnetic fields and the
associated concepts of energy, within the framework of Maxwell's equations.
Maxwell's equations mathematically explain how electric charges and currents
generate these fields and how changing fields induce each other. Formulas
for electromagnetic energy density and flux derived from these equations
demonstrate that electromagnetic fields carry energy and how this energy
propagates through space. Electromagnetic field theory is a cornerstone of
modern physics and engineering, with applications ranging from wireless
communications to optics. Much work has been done in this field, and we
present some of these studies.

In \cite{1} characterizes directional derivatives using an asymptotic
orthonormal frame and presents extended Serret-Frenet relations via cone
Frenet formulas. It explains the geometric meaning of energy on each
asymptotic orthonormal vector field in the null cone and expresses the
bending elastic energy for a particle based on its curve. The results are
supported by sketches showing energy variations with directional
derivatives. Additionally, it provides a geometric interpretation of energy
for unit vector fields and formulates Maxwell's equations for electric and
magnetic field vectors in null cone 3-space. Studies \cite{2,3} examine how
magnetic fields affect particle paths on a lightlike cone and characterize
magnetic curves using Killing magnetic fields. Studies \cite{4,10}
investigate the energy and volume of vector fields. In \cite{7}, this study
examines Berry's phase and defines Rytov parallel transport for
electromagnetic curves in an optic fiber using an alternative moving frame.
It also analyses electromagnetic curves with anholonomic coordinates for
Maxwellian evolution via Maxwell's equations. Studies \cite{8,12} explain
how geometric phase rotation relates to topological features in classical
Maxwell theory, using differential geometry to analyse various fiber paths.
In \cite{13}, the author explores the link between solutions of the cubic
non-linear Schrodinger equation and the localized induction equation. In 
\cite{14}, this study investigates the geometric properties of singular
Bertrand and Mannheim curves in 3D space forms. It also establishes
relationships between the singularities of these curves and the torsion of
their corresponding mate curves. Study \cite{15} describes a particle's
motion and calculates its bending elastic energy in 3D De-Sitter space. In 
\cite{16}, the authors examine the connection between electromagnetic theory
and Maxwell's equations. Study \cite{17} analyses how the Willmore energy of
curves in 3D Lorentzian space changes, describing variations in the Frenet
frame, curvature, and torsion. In study \cite{20}, the rotation of light's
polarization in a single-mode optical fiber following a curved path is
described. The study includes measurements of this rotation in a helical
fiber (bent into a spiral shape) with constant twist. In study \cite{21},
the authors present a geometric generalization of the action for a moving
particle's path in various spacetimes. In \cite{22}, the authors reduce a
hydrodynamics problem to a Heisenberg spin equation with constraints. In 
\cite{23}, the author examines how a magnetic field, generated by electric
current in an optical fiber current transformer, causes light polarization
to rotate as it travels through the fiber wrapped around a conductor. In 
\cite{25} the author states that the energy of a unit vector field on a
Riemannian manifold equals the energy of a related mapping on the unit
tangent bundle.

\section{Preliminaries}

Now we introduce some basic notions in semi-Euclidean space and curves. Let $%
\mathbb{R}
_{v}^{n+1}$ denote the $\left( n+1\right) -$dimensional pseudo-Euclidean
space of index $v\geqslant 0$; let $E=\left\{
e_{1},e_{2},...,e_{n+1}\right\} $ be an canonical \ basis of $%
\mathbb{R}
_{v}^{n+1}$. We choose two vectors $\varkappa ,\varrho \in 
\mathbb{R}
_{v}^{n+1}$, and the standard metric of $%
\mathbb{R}
_{v}^{n+1}$ is given by%
\begin{equation}
\left\langle \varkappa ,\varrho \right\rangle =-\overset{v}{\underset{i=1}{%
\sum }}\varkappa _{i}\varrho _{i}+\overset{n+1}{\underset{j=v+1}{\sum }}%
\varkappa _{j}\varrho _{j},  \tag{2.1}
\end{equation}%
where $\varkappa _{i}$ and $\varrho _{i}$ stand for the coordinate
components of $\varkappa $ and $\varrho $ with respect to $E$ in $%
\mathbb{R}
_{v}^{n+1}$, respectively.

For the vector $\varkappa \in 
\mathbb{R}
_{v}^{n+1}$, the vector $x$ is said to be spacelike if $\left\langle
\varkappa ,\varkappa \right\rangle >0$ or $\varkappa =0$, timelike if $%
\left\langle \varkappa ,\varkappa \right\rangle <0,$ lightlike(null) if $%
\left\langle \varkappa ,\varkappa \right\rangle =0$, $\varkappa \neq 0$.

Define the norm of a non-null vector $\varkappa $ by $\left\Vert \varkappa
\right\Vert =\left\vert \left\langle \varkappa ,\varkappa \right\rangle
\right\vert ^{\frac{1}{2}}$, where $\varkappa \in 
\mathbb{R}
_{v}^{n+1}$. We call $\varkappa $ the unit vector if $\left\Vert \varkappa
\right\Vert =1$.

Let $M_{q}^{3}$ $(c)\subset 
\mathbb{R}
_{v}^{n+1}$ denote the non flat 3D space forms of $q=0,1$ and constant
curvature $c\neq 0$. Meanwhile, $v=q$ if $c=1$, and $v=q+1,$ if $c=-1$.
Moreover, we will denote $M_{q}^{3}$ $(c)$ by the pseudo-Euclidean
hypersphere $S_{q}^{3}(1)$ or the pseudo-Euclidean hyperbolic space $%
H_{q}^{3}(-1)$ according to $c=1$ or $c=-1$, respectively, where $%
S_{q}^{3}(1)$ is denoted by%
\begin{equation*}
S_{q}^{3}(1)=\left\{ \varkappa =\left( \varkappa _{1},...,\varkappa
_{4}\right) \in 
\mathbb{R}
_{q}^{4}\mid \left\langle \varkappa ,\varkappa \right\rangle =1\right\} 
\end{equation*}%
and the pseudo-Euclidean hyperbolic space of index $q\geq 0$ and curvature $%
c=-1$ is given by 
\begin{equation*}
H_{q}^{3}(-1)=\{\varkappa =\left( \varkappa _{1},...,\varkappa _{4}\right)
\in 
\mathbb{R}
_{q+1}^{4}\mid \left\langle \varkappa ,\varkappa \right\rangle =-1\}.
\end{equation*}

Let $\gamma :I\rightarrow 
\mathbb{R}
_{v}^{n+1}$ be a curve in $%
\mathbb{R}
_{v}^{n+1}$ and let $\gamma ^{\prime }$ be the velocity vector of $\gamma $,
where $I$ is an open interval of $%
\mathbb{R}
$. For any $s\in I$, the curve $\gamma $ is called timelike curve, spacelike
curve or lightlike (null) curve if, for each $\left\langle \gamma ^{\prime
},\gamma ^{\prime }\right\rangle <0$, $\left\langle \gamma ^{\prime },\gamma
^{\prime }\right\rangle >0$ or $\left\langle \gamma ^{\prime },\gamma
^{\prime }\right\rangle =0$ and $\gamma ^{\prime }\neq 0$, respectively. We
call $\gamma $ a non null curve if $\gamma $ is a timelike curve or a
spacelike curve.

The Frenet frame of a non null curve in $M_{q}^{3}$ $(c)$ is as follows. Let 
$\gamma :I\rightarrow M_{q}^{3}$ $(c),$ $q=0,1$ be a non-null curve immersed
in the 3D space $M_{q}^{3}$ $(c)$, where $I$ is an open interval. If $%
\left\Vert \gamma ^{\prime }\right\Vert $ $=1$ for some $s\in I$, the curve $%
\gamma $ is called a unit speed curve. Then, in this paper $\gamma $ is
parametrized by the arc length parameter $s$. Letting $\nabla $ be the
Levi-Civita connection of $%
\mathbb{R}
_{v}^{4}$, there exists the Frenet frame $\left\{ T,N,B\right\} $ along $%
\gamma $ and smooth functions $\kappa ,\tau $ in $M_{q}^{3}$ $(c)$ such that%
\begin{equation*}
\nabla _{T}T=-\varepsilon _{1}c\gamma +\varepsilon _{2}\kappa N
\end{equation*}%
\begin{equation}
\nabla _{T}N=-\varepsilon _{1}\kappa T+\varepsilon _{3}\tau B  \tag{2.2}
\end{equation}%
\begin{equation*}
\nabla _{T}B=-\varepsilon _{2}\tau N,
\end{equation*}%
where $\kappa $ and $\tau $ are called the curvature and torsion of $\gamma $%
, respectively. Considering $\left\langle T,T\right\rangle $ $=\varepsilon
_{1},$ $\left\langle N,N\right\rangle $ $=\varepsilon _{2},$ $\left\langle
B,B\right\rangle $ $=\varepsilon _{3}$, and we denote by $\left\{
\varepsilon _{1},\varepsilon _{2},\varepsilon _{3}\right\} $ the casual
characters of $\left\{ T,N,B\right\} $. When $\left\{ T,N,B\right\} $ are
spacelike, then $\varepsilon _{i}=1$, and otherwise, $\varepsilon _{i}=-1$,
where $i\in \{1,2,3\}$. It is well known that curvature and torsion are
invariant under the isometries of $M_{q}^{3}$ $(c)$. Three vector fields $%
T,N,B$ consisting of the Frenet frame of $\gamma $ are called the tangent,
principal normal and binormal vector fields, respectively.

A vector field $M$ on $M_{q}^{3}$ $(c)$ along $\gamma $ is said to be
parallel along $\gamma $ if \ $\nabla _{s}M=0$, where $\nabla _{s}$ denotes
the covariant derivative along $\gamma $. A vector $M_{\gamma (s)}$ at $%
\gamma (s)$ is called parallel displacement of vector $M_{\gamma (s)}$ at $%
\gamma (s)$ along $\gamma $. If its tangent vector field $\gamma ^{\prime
}(s)$ of curve $\gamma $ is parallel along $\gamma $, then the curve is
called geodesic. We can denote the exponential map at $w\in M_{q}^{3}$ $(c)$
by $exp_{w}$ and review the exponential map $exp_{w}:T_{w}M_{q}^{3}$ $%
(c)\rightarrow M_{q}^{3}$ $(c)$, at $w\in M_{q}^{3}$ $(c)$ which is defined
by $exp_{w}(v)=\varsigma _{v}(1)$, where $\varsigma _{v}:[0,\infty
]\rightarrow $ $M_{q}^{3}$ $(c)$ is the constant speed geodesic starting
from $w$ with the initial velocity $\varsigma _{v}^{\prime }(0)=v$. For any
point $\gamma (s)$ in the curve $\gamma $, the principal normal geodesic in $%
M_{q}^{3}$ $(c)$ starting at $\gamma $ is defined as the geodesic curve $%
\varsigma _{s}^{\gamma }(t)=exp_{\gamma (s)}(tN(s))=f_{1}(t)\gamma
(s)+f_{2}(t)N(s),t\in 
\mathbb{R}
$, where the functions $f$ and $g$ are given by 
\begin{eqnarray*}
f_{1}(t) &=&\cos t,\text{ }f_{2}(t)=\sin t,\text{ if }\varepsilon _{2}c=1, \\
f_{1}(t) &=&\cosh t,\text{ }f_{2}(t)=\sinh t,\text{ if }\varepsilon _{2}c=-1,
\end{eqnarray*}%
\cite{14,18,19,24,27}.

\begin{definition}
For two Riemannian manifolds $(M,\varrho )$ and $\left( N,h\right) $ the
energy of a differentiable map $f:(M,\varrho )\rightarrow \left( N,h\right) $
is given as 
\begin{equation}
energy(f)=\frac{1}{2}\int_{M}\underset{a=1}{\overset{n}{\sum }}%
h(df(e_{a}),df(e_{a}))v,  \tag{2.3}
\end{equation}%
where $\left\{ e_{a}\right\} $ is a local basis of the tangent space and $v$
is the canonical volume form in $M$ \cite{25}.
\end{definition}

\begin{definition}
Let $Q:T(T^{1}M)\rightarrow T^{1}M$ be the connection map. Then, the
following conditions satisfy

i) $\omega oQ=\omega od\omega$\ and $\omega oQ=\omega o\varpi$ where $%
\varpi:T(T^{1}M)\rightarrow T^{1}M$ is the tangent bundle projection;

ii) for $\varrho \in T_{x}M$ and a section $\xi :M\rightarrow T^{1}M$; we
have 
\begin{equation}
Q(d\xi (\varrho ))=D_{\varrho }\xi ,  \tag{2.4}
\end{equation}%
where $D$ is the Levi-Civita covariant derivative \cite{25}.
\end{definition}

\begin{definition}
For $\varsigma_{1},\varsigma_{2}\in T_{\xi}\left( T^{1}M\right) $,
Riemannian metric on $TM$ is defined as 
\begin{equation}
\varrho_{S}(\varsigma_{1},\varsigma_{2})=\varrho(d\omega\left( \varsigma
_{1}\right) ,d\omega\left( \varsigma_{2}\right) )+\varrho(Q\left(
\varsigma_{1}\right) ,Q\left( \varsigma_{2}\right) ).  \tag{2.5}
\end{equation}

Here, as known $\varrho _{S}$ is called the Sasaki metric that also makes
the projection $\omega :T^{1}M\rightarrow M$ a Riemannian submersion \cite%
{25}.
\end{definition}

\section{The representation of the extended Serret-Frenet relations in non
flat 3-dimensional space forms $M_{q}^{3}$ $(c)$}

In this section, the directional derivatives are expressed in accordance
with the frame $\left\{ T,N,M\right\} $ in $M_{q}^{3}$ $(c)$ and the
extended Serret-Frenet relations are given using Frenet formulas. The
curvature of vector lines in anholonomic coordinates involves an additional
"twist" or "torsion" resulting not only from the metric properties of space
(e.g., length and angle) but also from the anholonomic constraints
themselves. The concepts of anholonomic coefficients, torsion, and
anholonomic connections are fundamental tools for understanding the
geometric properties and curvatures of vector lines in such systems, a way
to \ geometrically the complexity and path dependence of the system's paths
in state space are expressed.

Assuming that $\gamma =\gamma (s,\xi ,\eta )$ is a space curve lying in $%
M_{q}^{3}$ $(c)$, where $s$ is the distance along the $s$-lines of the curve
in the tangential direction so that unit tangent vector of $s$-lines is
defined by $T=T(s,\xi ,\eta )=\partial _{s}\gamma $, $N$ is the distance
along $\xi $-lines of the curve in the normal direction so that unit tangent
vector of $\xi $-lines is defined by $N=N(s,\xi ,\eta )=\partial _{\xi }N$, $%
B$ is the distance along the $B$-lines of the curve in the binormal
direction so that unit tangent vector of $B$-lines is defined by $B=B(s,\xi
,\eta )=\partial _{\eta }B$.

Hence, we can express the extended Serret-Frenet relations in $M_{q}^{3}$ $%
(c)$. First of all, to find the extended Frenet relations let's think the
the gradient operator $\nabla $ given by 
\begin{equation}
\nabla =\overrightarrow{T}\frac{\partial }{\partial s}+\overrightarrow{N}%
\frac{\partial }{\partial \xi }+\overrightarrow{B}\frac{\partial }{\partial
\eta },  \tag{3.1}
\end{equation}%
the curl and the divergence operator acting on an arbitrary vector $T$ is
written respectively, as%
\begin{equation}
\ Div \overrightarrow{T}=\nabla \overrightarrow{T}=\overrightarrow{T}%
\frac{\partial T}{\partial s}+\overrightarrow{N}\frac{\partial T}{\partial
\xi }+\overrightarrow{B}\frac{\partial T}{\partial \eta },  \tag{3.2}
\end{equation}%
\begin{equation}
\ Curl\overrightarrow{T}=\nabla \times T=\overrightarrow{T}\times \frac{%
\partial T}{\partial s}+\overrightarrow{N}\times \frac{\partial T}{\partial
\xi }+\overrightarrow{B}\times \frac{\partial T}{\partial \eta }.  \tag{3.3}
\end{equation}

First of all, we must create the Serre-Frenet frame according to the
parameters in the direction of the vector fields. The directional
derivatives along these unit vectors are defined by 
\begin{equation}
\frac{\partial }{\partial s}=T\cdot \nabla ,\frac{\partial }{\partial \xi }%
=N\cdot \nabla ,\frac{\partial }{\partial \eta }=B\cdot \nabla .  \tag{3.4}
\end{equation}

The directional derivatives can be obtained the tangential, principal normal
and binormal directions to the streamlines, respectively. For the
directional derivatives of the vector fields $T,N,B$ with respect to $\xi ,$
we can calculate as follows, for $a_{1}^{i}\in C^{\infty },$ $i=1,2,3.$

a) For $\frac{\partial \overrightarrow{T}}{\partial \xi },$ we have 
\begin{equation}
\frac{\partial \overrightarrow{T}}{\partial \xi }=a_{1}^{1}\overrightarrow{T}%
+a_{2}^{1}\overrightarrow{N}+a_{3}^{1}\overrightarrow{B}\Rightarrow \frac{%
\partial \overrightarrow{T}}{\partial \xi }=\varepsilon _{2}\Gamma
_{TN}^{\xi }\overrightarrow{N}+\varepsilon _{3}\Gamma _{TB}^{\xi }%
\overrightarrow{B}.  \tag{3.5}
\end{equation}

b) For $\frac{\partial \overrightarrow{N}}{\partial \xi },$ we have 
\begin{equation}
\frac{\partial \overrightarrow{N}}{\partial \xi }=a_{1}^{2}\overrightarrow{T}%
+a_{2}^{2}\overrightarrow{N}+a_{3}^{2}\overrightarrow{B}\Rightarrow \frac{%
\partial \overrightarrow{N}}{\partial \xi }=-\varepsilon _{1}\Gamma
_{TN}^{\xi }\overrightarrow{T}+\varepsilon _{3}\Gamma _{NB}^{\xi }%
\overrightarrow{B}.  \tag{3.6}
\end{equation}

c) For $\frac{\partial \overrightarrow{B}}{\partial \xi },$ we have 
\begin{equation}
\frac{\partial \overrightarrow{B}}{\partial \xi }=a_{1}^{3}\overrightarrow{T}%
+a_{2}^{3}\overrightarrow{N}+a_{3}^{3}\overrightarrow{B}\Rightarrow \frac{%
\partial \overrightarrow{B}}{\partial \xi }=-\varepsilon _{1}\Gamma
_{TB}^{\xi }\overrightarrow{T}-\varepsilon _{2}\Gamma _{NB}^{\xi }%
\overrightarrow{N}.  \tag{3.7}
\end{equation}
\ 

In this context, the directional derivatives of the vector fields $T,N,B$
with respect to $\eta $ for the given curve $\alpha $ are written as follows:%
\begin{equation}
\frac{d}{d\xi }%
\begin{bmatrix}
T \\ 
N \\ 
B%
\end{bmatrix}%
=%
\begin{bmatrix}
0 & \varepsilon _{2}\Gamma _{TN}^{\xi } & \varepsilon _{3}\Gamma _{TB}^{\xi }
\\ 
-\varepsilon _{1}\Gamma _{TN}^{\xi } & 0 & \varepsilon _{3}\Gamma _{NB}^{\xi
} \\ 
-\varepsilon _{1}\Gamma _{TB}^{\xi } & -\varepsilon _{2}\Gamma _{NB}^{\xi }
& 0%
\end{bmatrix}%
\begin{bmatrix}
T \\ 
N \\ 
B%
\end{bmatrix}%
.  \tag{3.8}
\end{equation}

By performing similar operations, the following equations is obtained,
respectively
\begin{equation}
\frac{d}{d\eta }%
\begin{bmatrix}
T \\ 
N \\ 
B%
\end{bmatrix}%
=%
\begin{bmatrix}
0 & \varepsilon _{2}\Upsilon _{TN}^{\eta } & \varepsilon _{3}\Upsilon
_{TB}^{\eta } \\ 
-\varepsilon _{1}\Upsilon _{TN}^{\eta } & 0 & \varepsilon _{3}\Upsilon
_{NB}^{\eta } \\ 
-\varepsilon _{1}\Upsilon _{TB}^{\eta } & -\varepsilon _{2}\Upsilon
_{NB}^{\eta } & 0%
\end{bmatrix}%
\begin{bmatrix}
T \\ 
N \\ 
B%
\end{bmatrix}%
,  \tag{3.9}
\end{equation}%
where $-\Upsilon _{NB}^{\xi }=\ Div \overrightarrow{B}$ by our
assumptions. 

We will now try to express the functions in (3.8) and (3.9). In summary,
other geometric quantities are computed by the vector analysis formulae in
the following manner. One find the followings;

a) For $\ Div\overrightarrow{T},$ since $\nabla _{T}T=\frac{\partial 
\overrightarrow{T}}{\partial s}=-\varepsilon _{1}c\gamma +\varepsilon
_{2}\kappa N,$ we get 
\begin{equation}
\ Div \overrightarrow{T}=\nabla \overrightarrow{T}=\overrightarrow{T}%
\frac{\partial T}{\partial s}+\overrightarrow{N}\frac{\partial T}{\partial
\xi }+\overrightarrow{B}\frac{\partial T}{\partial \eta }=\varepsilon
_{2}\Gamma _{TN}^{\xi }+\varepsilon _{3}\Upsilon _{TB}^{\eta }.  \tag{3.10}
\end{equation}

b) For $\ Div \overrightarrow{N},$ since $\nabla _{T}N=\frac{\partial 
\overrightarrow{N}}{\partial s}=-\varepsilon _{1}\kappa T+\varepsilon
_{3}\tau B,$ we get%
\begin{equation}
\ Div \overrightarrow{N}=\nabla \overrightarrow{N}=\overrightarrow{T}%
\frac{\partial N}{\partial s}+\overrightarrow{N}\frac{\partial N}{\partial
\xi }+\overrightarrow{B}\frac{\partial N}{\partial \eta }=-\varepsilon
_{1}\kappa +\varepsilon _{3}\Upsilon _{NB}^{\eta }.  \tag{3.11}
\end{equation}

c) For $\ Div \overrightarrow{B},$ since $\nabla _{T}B=\frac{\partial 
\overrightarrow{B}}{\partial s}=-\varepsilon _{2}\tau N,$ we get%
\begin{equation}
\ Div \overrightarrow{B}=\nabla \overrightarrow{B}=\overrightarrow{T}%
\frac{\partial B}{\partial s}+\overrightarrow{N}\frac{\partial B}{\partial
\xi }+\overrightarrow{B}\frac{\partial B}{\partial \eta }=-\Gamma _{NB}^{\xi
}.  \tag{3.12}
\end{equation}

Thus, we obtain 
\begin{equation}
\Gamma _{NB}^{\xi }=-\ Div \overrightarrow{B}=\varepsilon _{1}\kappa +%
\ Div \overrightarrow{N}.  \tag{3.13}
\end{equation}

On the other hand, we also obtain

d) For $\ Curl \overrightarrow{T},$ since $\frac{\partial 
\overrightarrow{T}}{\partial s}=-\varepsilon _{1}c\gamma +\varepsilon
_{2}\kappa N$ and by using equation $\ Curl \overrightarrow{T}$, we get%
\begin{equation}
\ Curl \overrightarrow{T}=-\varepsilon _{1}c\overrightarrow{T}\times
\gamma +\varepsilon _{2}\varepsilon _{3}\kappa \overrightarrow{B}%
+\varepsilon _{1}\left( \varepsilon _{3}\Gamma _{TB}^{\xi }-\varepsilon
_{2}\Upsilon _{TN}^{\eta }\right) \overrightarrow{T},  \tag{3.14}
\end{equation}%
where $\ Curl \overrightarrow{T}\cdot \overrightarrow{T}=\varepsilon
_{3}\Gamma _{TB}^{\xi }-\varepsilon _{2}\Upsilon _{TN}^{\eta }.$

e) For $\ Curl \overrightarrow{N},$ since $\frac{\partial 
\overrightarrow{N}}{\partial s}=-\varepsilon _{1}\kappa T+\varepsilon
_{3}\tau B$ and for the equation $\ Curl \overrightarrow{N}$, we get 
\begin{equation}
\ Curl \overrightarrow{N}=\varepsilon _{1}\varepsilon _{3}\Gamma
_{NB}^{\xi }\overrightarrow{T}-\varepsilon _{2}\left( \varepsilon _{3}\tau
+\varepsilon _{1}\Upsilon _{TN}^{\eta }\right) \overrightarrow{N}%
+\varepsilon _{1}\varepsilon _{3}\Gamma _{TN}^{\xi }\overrightarrow{B}, 
\tag{3.15}
\end{equation}%
where $\ Curl \overrightarrow{N}\cdot \overrightarrow{N}=-\varepsilon
_{3}\tau -\varepsilon _{1}\Upsilon _{TN}^{\eta }.$

g) For $\ Curl \overrightarrow{B},$ since $\frac{\partial 
\overrightarrow{B}}{\partial s}=-\varepsilon _{2}\tau N,$for the equation $%
\ Curl \overrightarrow{B},$we have 
\begin{equation}
\ Curl \overrightarrow{B}=\varepsilon _{1}\varepsilon _{2}\Upsilon
_{NB}^{\eta }\overrightarrow{T}-\varepsilon _{2}\left( \varepsilon _{2}\tau
+\varepsilon _{1}\Upsilon _{TB}^{\eta }\right) \overrightarrow{N}%
+\varepsilon _{1}\varepsilon _{3}\Gamma _{TB}^{\xi }\overrightarrow{B}, 
\tag{3.16}
\end{equation}%
where $\ Curl \overrightarrow{B}\cdot \overrightarrow{B}=\varepsilon
_{1}\Gamma _{TB}^{\xi }.$ Therefore, we get 
\begin{equation}
\frac{\partial }{\partial \xi }\overrightarrow{T}\cdot \overrightarrow{N}%
=\Gamma _{TN}^{\xi };\frac{\partial }{\partial \xi }\overrightarrow{T}\cdot 
\overrightarrow{B}=\Gamma _{TB}^{\xi };\frac{\partial }{\partial \xi }%
\overrightarrow{N}\cdot \overrightarrow{B}=\Gamma _{NB}^{\xi };\ Div%
\overrightarrow{B}=-\Gamma _{NB}^{\xi };  \tag{3.17a}
\end{equation}%
\begin{equation}
\Psi _{B}=\ Curl \overrightarrow{B}\cdot \overrightarrow{B}=\varepsilon
_{1}\Gamma _{TB}^{\xi };\Psi _{N}=\ Curl\overrightarrow{N}\cdot 
\overrightarrow{N}=-\varepsilon _{3}\tau -\varepsilon _{1}\Upsilon
_{TN}^{\eta };  \tag{3.17b}
\end{equation}%
\begin{equation}
\Psi _{T}=\ Curl\overrightarrow{T}\cdot \overrightarrow{T}=\varepsilon
_{3}\Gamma _{TB}^{\xi }-\varepsilon _{2}\Upsilon _{TN}^{\eta }  \tag{3.17c}
\end{equation}%
and some functions can be given as 
\begin{equation}
\Gamma _{TB}^{\xi }=\varepsilon _{1}\ Curl\overrightarrow{B}\cdot 
\overrightarrow{B};\Upsilon _{TN}^{\eta }=-\varepsilon _{1}\varepsilon
_{3}\tau -\varepsilon _{1}\ Curl\overrightarrow{N}\cdot \overrightarrow{%
N}.  \tag{3.18}
\end{equation}

This implies%
\begin{equation}
\ Curl\overrightarrow{N}\cdot \overrightarrow{B}=\varepsilon _{1}\Gamma
_{TN}^{\xi };\ Curl \overrightarrow{N}\cdot \overrightarrow{T}%
=\varepsilon _{3}\Gamma _{NB}^{\xi }=-\varepsilon _{3}\ Div%
\overrightarrow{B},  \tag{3.19a}
\end{equation}%
\begin{equation}
\ Curl \overrightarrow{B}\cdot \overrightarrow{N}=-\varepsilon _{2}\tau
-\varepsilon _{1}\Upsilon _{TB}^{\eta };\ Curl\overrightarrow{B}\cdot 
\overrightarrow{T}=\varepsilon _{2}\Upsilon _{NB}^{\eta };\ Curl%
\overrightarrow{T}\cdot \overrightarrow{B}=\varepsilon _{2}\kappa , 
\tag{3.19b}
\end{equation}%
\begin{equation}
\Upsilon _{NB}^{\eta }=\varepsilon _{2}\ Curl\overrightarrow{B}\cdot 
\overrightarrow{T};\Upsilon _{TB}^{\eta }=-\varepsilon _{1}\ Curl%
\overrightarrow{B}\cdot \overrightarrow{N}-\varepsilon _{1}\varepsilon
_{2}\tau ;\Upsilon _{TN}^{\xi }=\varepsilon _{1}\ Curl\overrightarrow{N}%
\cdot \overrightarrow{B}.  \tag{3.20}
\end{equation}

Therefore, from the last equations, if we substitute the obtained values of
the smooth functions, we write Serret-Frenet relations in the following forms%
\begin{equation}
\frac{d}{d\xi }%
\begin{bmatrix}
T \\ 
N \\ 
B%
\end{bmatrix}%
=%
\begin{bmatrix}
0 & -\varepsilon _{1}\varepsilon _{3}\ Curl\overrightarrow{N}\cdot 
\overrightarrow{B} & \varepsilon _{1}\varepsilon _{3}\Psi _{B} \\ 
-\ Curl\overrightarrow{N}\cdot \overrightarrow{B} & 0 & -\varepsilon
_{3}\ Div\overrightarrow{B} \\ 
-\Psi _{B} & -\varepsilon _{2}\ Div\overrightarrow{B} & 0%
\end{bmatrix}%
\begin{bmatrix}
T \\ 
N \\ 
B%
\end{bmatrix}
\tag{3.21}
\end{equation}%
and%
\begin{equation}
\frac{d}{d\eta }%
\begin{bmatrix}
T \\ 
N \\ 
B%
\end{bmatrix}%
=%
\begin{bmatrix}
0 & -\varepsilon _{1}\varepsilon _{3}\tau -\varepsilon _{1}\Psi _{N} & 
\begin{array}{c}
-\varepsilon _{1}\varepsilon _{3}(\varepsilon _{2}\tau  \\ 
+\ Curl\overrightarrow{B}\cdot \overrightarrow{N})%
\end{array}
\\ 
\varepsilon _{3}\tau +\Psi _{N} & 0 & \varepsilon _{2}\varepsilon _{3}\ Curl\overrightarrow{B}\cdot \overrightarrow{T} \\ 
\begin{array}{c}
\varepsilon _{2}\tau  \\ 
+\ Curl\overrightarrow{B}\cdot \overrightarrow{N}%
\end{array}
& -\ Curl\overrightarrow{B}\cdot \overrightarrow{T} & 0%
\end{bmatrix}%
\begin{bmatrix}
T \\ 
N \\ 
B%
\end{bmatrix}%
,  \tag{3.22}
\end{equation}%
where $\kappa $ is the curvature function and $\tau $ is the torsion
function of the unit speed timelike curve $\gamma (s,\xi ,\eta )$.

This relations was originally obtained and application of the identity $%
\ Curl\nabla h=0$ yields
\begin{eqnarray*}
\ Curl\nabla h &=&\overrightarrow{T}\times \frac{\partial \nabla h}{%
\partial s}+\overrightarrow{N}\times \frac{\partial \nabla h}{\partial \xi }+%
\overrightarrow{B}\times \frac{\partial \nabla h}{\partial \eta } \\
&=&\frac{\partial h}{\partial s}\ Curl\overrightarrow{T}+\frac{\partial
h}{\partial \xi }\ Curl\overrightarrow{N}+\frac{\partial h}{\partial
\eta }\ Curl\overrightarrow{B}+\overrightarrow{T}\times (%
\overrightarrow{T}\frac{\partial ^{2}h}{\partial s^{2}}+\overrightarrow{N}%
\frac{\partial ^{2}h}{\partial s\partial \xi }+\overrightarrow{B}\frac{%
\partial ^{2}h}{\partial s\partial \eta })
\end{eqnarray*}%
\begin{equation*}
+\overrightarrow{N}\times (\overrightarrow{T}\frac{\partial ^{2}h}{\partial
\xi \partial s}+\overrightarrow{N}\frac{\partial ^{2}h}{\partial \xi ^{2}}+%
\overrightarrow{B}\frac{\partial ^{2}h}{\partial \xi \partial \eta })+%
\overrightarrow{B}\times (\overrightarrow{T}\frac{\partial ^{2}h}{\partial
\eta \partial s}+\overrightarrow{N}\frac{\partial ^{2}h}{\partial \eta
\partial \xi }+\overrightarrow{B}\frac{\partial ^{2}h}{\partial \eta ^{2}})
\end{equation*}%
\begin{eqnarray*}
&=&\frac{\partial h}{\partial s}\ Curl\overrightarrow{T}+\frac{\partial
h}{\partial \xi }\ Curl\overrightarrow{N}+\frac{\partial h}{\partial
\eta }\ Curl\overrightarrow{B}+\varepsilon _{3}\left( \frac{\partial
^{2}h}{\partial s\partial \xi }-\frac{\partial ^{2}h}{\partial \xi \partial s%
}\right) \overrightarrow{B} \\
&&+\varepsilon _{2}\left( \frac{\partial ^{2}h}{\partial \eta \partial s}-%
\frac{\partial ^{2}h}{\partial s\partial \eta }\right) \overrightarrow{N}%
+\varepsilon _{1}\left( \frac{\partial ^{2}h}{\partial \xi \partial \eta }-%
\frac{\partial ^{2}h}{\partial \eta \partial \xi }\right) \overrightarrow{T},
\end{eqnarray*}%
from the equations (3.14), (3.15), (3.16) and considering the property $%
\frac{\partial ^{2}}{\partial \xi \partial s}=\frac{\partial ^{2}}{\partial
s\partial \xi }$ for any two different parameters, we can write as follows%
\begin{equation*}
0=\frac{\partial h}{\partial s}\left( 
\begin{array}{c}
-\varepsilon _{1}c\overrightarrow{T}\times \gamma +\varepsilon
_{2}\varepsilon _{3}\kappa \overrightarrow{B} \\ 
+\varepsilon _{1}(\varepsilon _{3}\Gamma _{TB}^{\xi }-\varepsilon
_{2}\Upsilon _{TN}^{\eta })\overrightarrow{T}%
\end{array}%
\right) +\frac{\partial h}{\partial \xi }\left( 
\begin{array}{c}
\varepsilon _{1}\varepsilon _{3}\Gamma _{NB}^{\xi }\overrightarrow{T}%
-\varepsilon _{2}\left( \varepsilon _{3}\tau +\varepsilon _{1}\Upsilon
_{TN}^{\eta }\right) \overrightarrow{N} \\ 
+\varepsilon _{1}\varepsilon _{3}\Gamma _{TN}^{\xi }\overrightarrow{B}%
\end{array}%
\right) 
\end{equation*}%
\begin{equation}
+\frac{\partial h}{\partial \eta }\left( \varepsilon _{1}\varepsilon
_{2}\Upsilon _{NB}^{\eta }\overrightarrow{T}-\varepsilon _{2}\left(
\varepsilon _{2}\tau +\varepsilon _{1}\Upsilon _{TB}^{\eta }\right) 
\overrightarrow{N}+\varepsilon _{1}\varepsilon _{3}\Gamma _{TB}^{\xi }%
\overrightarrow{B}\right)   \tag{3.23}
\end{equation}%
and 
\begin{eqnarray*}
0 &=&-\varepsilon _{1}\frac{\partial h}{\partial s}\left( c\overrightarrow{T}%
\times \gamma \right)  \\
&&+\varepsilon _{1}\left( \frac{\partial ^{2}h}{\partial \xi \partial \eta }-%
\frac{\partial ^{2}h}{\partial \eta \partial \xi }+\frac{\partial h}{%
\partial s}(\varepsilon _{3}\Gamma _{TB}^{\xi }-\varepsilon _{2}\Upsilon
_{TN}^{\eta })+\frac{\partial h}{\partial \xi }\varepsilon _{3}\Gamma
_{NB}^{\xi }+\frac{\partial h}{\partial \eta }\varepsilon _{2}\Upsilon
_{NB}^{\eta }\right) \overrightarrow{T} \\
&&+\varepsilon _{2}\left( \frac{\partial ^{2}h}{\partial \eta \partial s}-%
\frac{\partial ^{2}h}{\partial s\partial \eta }+\frac{\partial h}{\partial
\xi }\left( -\left( \varepsilon _{3}\tau +\varepsilon _{1}\Upsilon
_{TN}^{\eta }\right) \right) +\frac{\partial h}{\partial \eta }\left(
-\varepsilon _{2}\tau -\varepsilon _{2}\varepsilon _{1}\Upsilon _{TB}^{\eta
}\right) \right) \overrightarrow{N} \\
&&+\varepsilon _{3}\left( \frac{\partial ^{2}h}{\partial s\partial \xi }-%
\frac{\partial ^{2}h}{\partial \xi \partial s}+\frac{\partial h}{\partial s}%
\varepsilon _{2}\kappa +\frac{\partial h}{\partial \xi }\varepsilon
_{1}\Gamma _{TN}^{\xi }+\frac{\partial h}{\partial \eta }\varepsilon
_{1}\Gamma _{TB}^{\xi }\right) \overrightarrow{B}.
\end{eqnarray*}

If the algebraic equality is taken into account from the last equations
above, the following equation system can be written%
\begin{equation}
\frac{\partial ^{2}h}{\partial \xi \partial s}-\frac{\partial ^{2}h}{%
\partial s\partial \xi }=\frac{\partial h}{\partial s}\varepsilon _{2}\kappa
+\frac{\partial h}{\partial \xi }\varepsilon _{1}\ Curl \overrightarrow{N%
}\cdot \overrightarrow{B}+\frac{\partial h}{\partial \eta }\varepsilon _{1}%
\ Curl \overrightarrow{B}\cdot \overrightarrow{B}  \tag{3.24a}
\end{equation}%
\begin{equation}
\frac{\partial ^{2}h}{\partial s\partial \eta }-\frac{\partial ^{2}h}{%
\partial \eta \partial s}=\frac{\partial h}{\partial \xi }\ Curl%
\overrightarrow{N}\cdot \overrightarrow{N}+\frac{\partial h}{\partial \eta }%
\ Curl\overrightarrow{B}\cdot \overrightarrow{N}  \tag{3.24b}
\end{equation}%
\begin{equation}
\frac{\partial ^{2}h}{\partial \eta \partial \xi }-\frac{\partial ^{2}h}{%
\partial \xi \partial \eta }=\frac{\partial h}{\partial s}\left( 
\begin{array}{c}
\varepsilon _{3}\ Curl \overrightarrow{B}\cdot \overrightarrow{B} \\ 
+\varepsilon _{2}(\varepsilon _{3}\tau  \\ 
+\varepsilon _{1}\ Curl \overrightarrow{N}\cdot \overrightarrow{N})%
\end{array}%
\right) -\varepsilon _{3}\frac{\partial h}{\partial \xi }\ Div%
\overrightarrow{B}+\frac{\partial h}{\partial \eta }\ Curl%
\overrightarrow{B}\cdot \overrightarrow{T}.  \tag{3.24c}
\end{equation}

Thus, considering the equations in (3.24), the following equations can be
written%
\begin{eqnarray*}
\frac{\partial h}{\partial s}\varepsilon _{2}\kappa +\frac{\partial h}{%
\partial \xi }\varepsilon _{1}\ Curl \overrightarrow{N}\cdot 
\overrightarrow{B}+\frac{\partial h}{\partial \eta }\varepsilon _{1}\ Curl\overrightarrow{B}\cdot \overrightarrow{B} &=&0 \\
\frac{\partial h}{\partial \xi }\ Curl \overrightarrow{N}\cdot 
\overrightarrow{N}+\frac{\partial h}{\partial \eta }\ Curl%
\overrightarrow{B}\cdot \overrightarrow{N} &=&0 \\
\frac{\partial h}{\partial s}(\varepsilon _{3}\ Curl \overrightarrow{B}%
\cdot \overrightarrow{B}+\varepsilon _{2}(\varepsilon _{3}\tau +\varepsilon
_{1}\ Curl \overrightarrow{N}\cdot \overrightarrow{N})) &=&\varepsilon
_{3}\frac{\partial h}{\partial \xi }\ Div \overrightarrow{B}-\frac{%
\partial h}{\partial \eta }\ Curl \overrightarrow{B}\cdot 
\overrightarrow{T}.
\end{eqnarray*}

\section{The Maxwell's equations of electromagnetic wave vector fields in $%
M_{q}^{3}\left( c\right) $}

It states that the orientation of an electromagnetic wave within an optical
fiber is defined using an orthogonal unit vector frame consisting of the
vector fields $\overrightarrow{T}$, $\overrightarrow{N}$ and $%
\overrightarrow{B}$. Orientation of the electromagnetic wave: This refers to
the properties of the electromagnetic wave in space $M_{q}^{3}\left(
c\right) $, such as its position, direction, or polarization state. The
directions of the wave's electric and magnetic field vectors as it travels
through the fiber are important. As an electromagnetic wave propagates
through an optical fiber, a geometric phase called the Berry phase arises
when the wave's vector fields or related parameters in specific $\xi $ and $%
\eta $ directions (possibly within the fiber's cross-section or related to
its polarization) change. This implies that the wave's motion within the
fiber not only acquires a dynamic phase but also acquires an additional
"geometric memory" as a result of the wave's spatial or polarization
structure following specific paths. Also, the Berry phase, a path-dependent
phase phenomenon associated with electromagnetic waves in optical fiber.
This implies that the phase is related to the wave's behaviour in specific
directions within the fiber's cross-section. This phase is a special type of
phase that occurs during the evolution of a quantum system or
electromagnetic waves). Normally, the phase change is related to the
system's energy and time (dynamic phase). However, the Berry phase depends
on the path followed by the system in parameter space. This phase depends on
the "geometry" of the path (the area it encloses in parameter space), not
the time itself or the energy. This is why it is called the "geometric"
phase.

In the optical context, this can occur when parameters such as the
polarization or orientation of light are slowly changed. These terms may
refer to directions defined in a specific context. Since the sentence refers
to the propagation of an electromagnetic wave along an optical fiber, these
directions may relate to the fiber's cross-section or the wave's
polarization. Generally, these terms may refer to parameters associated with
the components of the electromagnetic wave's vector field (electric or
magnetic field) in different directions within the fiber's cross-section or
its polarization state. Berry phase occurs when the wave's parameters change
along these directions within the fiber cross-section.

It states that the electric and magnetic field vectors ($\overrightarrow{E}$
and $\overrightarrow{M}$) of an electromagnetic wave propagating in an
optical fiber exhibit a rotation along the fiber axis (in the tangential $s-$%
direction) with respect to the $\{T,N,B\}$ reference frame defined by the
geometry of the fiber. This rotation can be caused by bending, torsion, or
polarization-related effects of the wave. This is an important phenomenon
for understanding polarization preservation or change in optical fiber.

Optical fiber can be defined as a curve $\gamma (s,\xi ,\eta )$ via
alternative moving frame in three dimensional space. If we want to
understand the electromagnetic theory, we have to know Maxwell's equations.
So that Electromagnetic waves propagated along the optical fiber and the
electromagnetic waves spread through the optical fiber in which its axis is
expressed by the curve $\gamma $. On account of the vectorial nature of the
light electromagnetic waves are defined by using the vector fields. The
orientation of the electromagnetic wave in the fiber is defined by using the
frame of vectors $\{T,N,B\}$ in $M_{q}^{3}\left( c\right) $.

For an electromagnetic wave of a space curve $\gamma $, the electric field
vector $\overrightarrow{E}$ and the magnetic field vector $\overrightarrow{M}
$ are expected to perform a rotation in the tangential direction according
to the unit vectors $\left\{ T,N,B\right\} $. Also, the electromagnetic wave
carries magnetic vector field $\overrightarrow{M}$. Consequently, the
electromagnetic vectors $\overrightarrow{E}$ and $\overrightarrow{M}$ may be
considered as a physically coordinate frame, which are expressed according
to orthonormal unit vectors $\left\{ T,N,B\right\} $.

We know that Maxwell's equations are a set of four partial differential
equations that form the basis of classical electromagnetism. They describe
how electric and magnetic fields behave and interact with each other and
with charges and currents. These equations demonstrate that light consists
of electromagnetic waves. Thus, the following equations are given for the
magnetic vector fields and the electric vector fields in our study.
\begin{equation}
\nabla \overrightarrow{E}_{s\xi \eta }^{\xi }=0;\nabla \overrightarrow{E}%
_{s\xi \eta }^{\eta }=0;\nabla \overrightarrow{M}^{\xi }=0;\nabla 
\overrightarrow{M}^{\eta }=0.  \tag{4.1}
\end{equation}

Let $\overrightarrow{E}$ and $\overrightarrow{M}$ be the vectors of the
electromagnetic wave, so that $\overrightarrow{E}$ and $\overrightarrow{M}$
are perpendicular to the tangent vector field $T=\gamma ^{\prime }$ along
the curve $\gamma \left( s,\xi ,\eta \right) $ \cite{9}.

We consider the fundamental fiber mode in the $\xi -$direction along the
optical fiber $\gamma $ according to frame $\left\{ T,N,B\right\} $ in $%
M_{q}^{3}\left( c\right) ,$ then 
\begin{equation*}
\left\langle \overrightarrow{E}^{\xi },\overrightarrow{T}\right\rangle =0.
\end{equation*}

The derivation of the electric vector $\overrightarrow{E}^{\xi }$ between
any two points in the $\xi -$direction along optical fiber $\gamma $  with
respect to frame $\left\{ T,N,B\right\} $ is given as
\begin{equation}
\frac{\partial \overrightarrow{E}^{\xi }}{\partial \xi }=c_{1}^{\xi }%
\overrightarrow{T}+c_{2}^{\xi }\overrightarrow{N}+c_{3}^{\xi }%
\overrightarrow{B},  \tag{4.2}
\end{equation}%
where $c_{i}^{\xi },$ $i=1,2,3$ are smooth functions.

The electric field vector$\overrightarrow{E}^{\xi }$ is perpendicular to the
vector $\overrightarrow{T}$ in the frame $\{T,N,B\}$, the vector $%
\overrightarrow{T}$ is tangent to the fiber axis or the wave's direction of
propagation. This means that the electric field vector $\overrightarrow{E}%
^{\xi }$ is at a 90 degree angle to this direction $\overrightarrow{T}$.
Recalling that electromagnetic waves are transverse waves in free space ($%
\overrightarrow{E}^{\xi }$ and $\overrightarrow{M}^{_{\xi }}$ are
perpendicular to the direction of propagation), this statement indicates
that the wave retains its transverse character within the fiber or that a
particular mode is transversely polarized. Therefore, since $\overrightarrow{%
E}^{\xi }$ and $\overrightarrow{M}^{_{\xi }}$ are perpendicular to the
tangent vector field $T=\gamma ^{\prime }$ along $\gamma \left( s,\xi ,\eta
\right) $, we have
\begin{equation}
\overrightarrow{T}\cdot \overrightarrow{E}^{\xi }=0,\overrightarrow{E}^{\xi
}\cdot \overrightarrow{E}^{\xi }=const\text{.};\overrightarrow{T}\cdot \frac{%
\partial \overrightarrow{E}^{\xi }}{\partial \xi }=-\overrightarrow{E}^{\xi
}\cdot \frac{\partial \overrightarrow{T}}{\partial \xi },\overrightarrow{E}%
^{\xi }\cdot \frac{\partial \overrightarrow{E}^{\xi }}{\partial \xi }=0 
\tag{4.3}
\end{equation}%
\begin{equation}
\overrightarrow{T}\cdot \overrightarrow{M}^{_{\xi }}=0,\overrightarrow{M}%
^{_{\xi }}\cdot \overrightarrow{M}^{_{\xi }}=const\text{.};\overrightarrow{T}%
\cdot \frac{\partial \overrightarrow{M}^{_{\xi }}}{\partial \xi }=-%
\overrightarrow{M}^{_{\xi }}\cdot \frac{\partial \overrightarrow{T}}{%
\partial \xi },\overrightarrow{M}^{_{\xi }}\cdot \frac{\partial 
\overrightarrow{M}^{_{\xi }}}{\partial \xi }=0  \tag{4.4}
\end{equation}%
and considering (4.3), we write $\overrightarrow{E}^{\xi }=E_{\xi }^{1}%
\overrightarrow{N}+E_{\xi }^{3}\overrightarrow{B},$ for the components of
the electric vector field and by using the equation (3.8)(or (3.21)), we
obtain%
\begin{eqnarray*}
c_{1}^{\xi } &=&\overrightarrow{T}\cdot \frac{\partial \overrightarrow{E}%
^{\xi }}{\partial \xi }\varepsilon _{1}=-\overrightarrow{E}^{\xi }\cdot 
\frac{\partial \overrightarrow{T}}{\partial \xi }\varepsilon
_{1}=-\varepsilon _{1}\left( E_{\xi }^{1}\Gamma _{TN}^{\xi }+E_{\xi
}^{3}\Gamma _{TB}^{\xi }\right)  \\
c_{2}^{\xi } &=&\overrightarrow{N}\cdot \frac{\partial \overrightarrow{E}%
^{\xi }}{\partial \xi }\varepsilon _{2}=-\overrightarrow{E}^{\xi }\cdot 
\frac{\partial \overrightarrow{N}}{\partial \xi }\varepsilon
_{2}=-\varepsilon _{2}E_{\xi }^{3}\Gamma _{NB}^{\xi } \\
c_{3}^{\xi } &=&\overrightarrow{B}\cdot \frac{\partial \overrightarrow{E}%
^{\xi }}{\partial \xi }\varepsilon _{3}=-\overrightarrow{E}^{\xi }\cdot 
\frac{\partial \overrightarrow{B}}{\partial \xi }\varepsilon
_{3}=\varepsilon _{3}E_{\xi }^{1}\Gamma _{NB}^{\xi }.
\end{eqnarray*}

Thus, if the values in the previous equations are taken into account in the
equation (4.2), we obtain 
\begin{equation}
\frac{\partial \overrightarrow{E}_{s\xi \eta }^{\xi }}{\partial \xi }%
=-\varepsilon _{1}\left( E_{\xi }^{1}\Gamma _{TN}^{\xi }+E_{\xi }^{3}\Gamma
_{TB}^{\xi }\right) \overrightarrow{T}-\varepsilon _{2}E_{\xi }^{3}\Gamma
_{NB}^{\xi }\overrightarrow{N}+\varepsilon _{3}E_{\xi }^{1}\Gamma _{NB}^{\xi
}\overrightarrow{B}.  \tag{4.5}
\end{equation}

The change of the electric vector field $\overrightarrow{E}^{\eta }$ with
respect to $\eta -$direction \ $\frac{\partial \overrightarrow{E}^{\eta }}{%
\partial \eta }$, we can write%
\begin{equation}
\frac{\partial \overrightarrow{E}^{\eta }}{\partial \eta }=c_{1}^{\eta }%
\overrightarrow{T}+c_{2}^{\eta }\overrightarrow{N}+c_{3}^{\eta }%
\overrightarrow{B},  \tag{4.6}
\end{equation}%
where $c_{i}^{\eta },$ $i=1,2,3$ are smooth functions. Also, the following
equations hold 
\begin{equation}
\overrightarrow{T}\cdot \overrightarrow{E}^{\eta }=0,\overrightarrow{E}%
^{\eta }\cdot \overrightarrow{E}^{\eta }=const;\overrightarrow{T}\cdot \frac{%
\partial \overrightarrow{E}^{\eta }}{\partial \eta }=-\overrightarrow{E}%
^{\eta }\cdot \frac{\partial \overrightarrow{T}}{\partial \eta },  \tag{4.7a}
\end{equation}%
\begin{equation}
\overrightarrow{E}^{\eta }\cdot \frac{\partial \overrightarrow{E}^{\eta }}{%
\partial \eta }=0,\text{ }\overrightarrow{E}^{\eta }=E_{\eta }^{1}%
\overrightarrow{N}+E_{\eta }^{3}\overrightarrow{B}.  \tag{4.7b}
\end{equation}

Hence, from the derivatives of the vector fields (3.9)(or (3.22)), we get 
\begin{equation}
c_{1}^{\eta }=-\varepsilon _{1}\left( E_{\eta }^{1}\Upsilon _{TN}^{\eta
}+E_{\eta }^{1}\Upsilon _{TB}^{\eta }\right) ;\text{ }c_{2}^{\eta
}=-\varepsilon _{2}E_{\eta }^{3}\Upsilon _{NB}^{\eta };\text{ }c_{3}^{\eta
}=\varepsilon _{3}E_{\eta }^{1}\Upsilon _{NB}^{\eta }  \tag{4.8a}
\end{equation}%
from the equations (4.8), we obtain 
\begin{equation}
\frac{\partial \overrightarrow{E}^{\eta }}{\partial \eta }=-\varepsilon
_{1}\left( E_{\eta }^{1}\Upsilon _{TN}^{\eta }+E_{\eta }^{1}\Upsilon
_{TB}^{\eta }\right) \overrightarrow{T}-\varepsilon _{2}E_{\eta
}^{3}\Upsilon _{NB}^{\eta }\overrightarrow{N}+\varepsilon _{3}E_{\eta
}^{1}\Upsilon _{NB}^{\eta }\overrightarrow{B}.  \tag{4.9}
\end{equation}

Similarly, for the change of the electric vector field $\overrightarrow{E}%
^{s}$ with respect to $s-$direction $\frac{\partial \overrightarrow{E}^{s}}{%
\partial s}$, we obtain 
\begin{equation}
\frac{\partial \overrightarrow{E}^{s}}{\partial s}=c_{1}^{s}\overrightarrow{T%
}+c_{2}^{s}\overrightarrow{N}+c_{3}^{s}\overrightarrow{B}  \tag{4.10}
\end{equation}%
\begin{equation}
\overrightarrow{E}^{s}\cdot \frac{\partial \overrightarrow{E}^{s}}{\partial s%
}=c_{2}^{s}\overrightarrow{E}^{s}.\overrightarrow{N}+c_{3}^{s}%
\overrightarrow{E}^{s}.\overrightarrow{B}=0\text{ and }\overrightarrow{E}%
^{s}=E_{s}^{1}\overrightarrow{N}+E_{s}^{3}\overrightarrow{B},  \tag{4.11}
\end{equation}%
where $c_{i}^{s},$ $i=1,2,3$ are smooth functions.

Therefore, from (2.2) the components of $\frac{\partial \overrightarrow{E}%
^{s}}{\partial s}$ are obtained as follows 
\begin{equation*}
c_{1}^{s}=-\varepsilon _{1}\kappa E_{s}^{1};\text{ \ }c_{2}^{s}=-\varepsilon
_{2}E_{s}^{3}\tau ;\text{ \ }c_{3}^{s}=\varepsilon _{3}E_{s}^{1}\tau 
\end{equation*}%
and when the last equations obtained are used, the following equation is
obtained 
\begin{equation}
\frac{\partial \overrightarrow{E}^{s}}{\partial s}=-\varepsilon _{1}\kappa
E_{s}^{1}\overrightarrow{T}-\varepsilon _{2}E_{s}^{3}\tau \overrightarrow{N}%
+\varepsilon _{3}E_{s}^{1}\tau \overrightarrow{B}.  \tag{4.12}
\end{equation}

Hence, we compute that 
\begin{align*}
\nabla \overrightarrow{E}& =\overrightarrow{T}\cdot \frac{\partial 
\overrightarrow{E}^{s}}{\partial s}+\overrightarrow{N}\cdot \frac{\partial 
\overrightarrow{E}^{\xi }}{\partial \xi }+\overrightarrow{B}\cdot \frac{%
\partial \overrightarrow{E}^{\eta }}{\partial \eta } \\
& =\overrightarrow{T}\cdot (-\varepsilon _{1}\kappa E_{s}^{1}\overrightarrow{%
T}-\varepsilon _{2}E_{s}^{3}\tau \overrightarrow{N}+\varepsilon
_{3}E_{s}^{1}\tau \overrightarrow{B}) \\
& +\overrightarrow{N}\cdot \left( -\varepsilon _{1}(E_{\xi }^{1}\Gamma
_{TN}^{\xi }+E_{\xi }^{3}\Gamma _{TB}^{\xi })\overrightarrow{T}-\varepsilon
_{2}E_{\xi }^{3}\Gamma _{NB}^{\xi }\overrightarrow{N}+\varepsilon _{3}E_{\xi
}^{1}\Gamma _{NB}^{\xi }\overrightarrow{B}\right)  \\
& +\overrightarrow{B}\cdot \left( -\varepsilon _{1}\left( E_{\eta
}^{1}\Upsilon _{TN}^{\eta }+E_{\eta }^{1}\Upsilon _{TB}^{\eta }\right) 
\overrightarrow{T}-\varepsilon _{2}E_{\eta }^{3}\Upsilon _{NB}^{\eta }%
\overrightarrow{N}+\varepsilon _{3}E_{\eta }^{1}\Upsilon _{NB}^{\eta }%
\overrightarrow{B}\right) 
\end{align*}%
\begin{equation}
\nabla \overrightarrow{E}=-\kappa E_{s}^{1}+E_{\eta }^{1}\Upsilon
_{NB}^{\eta }-E_{\xi }^{3}\Gamma _{NB}^{\xi }  \tag{4.13}
\end{equation}%
which implies that 
\begin{equation}
\nabla \overrightarrow{E}=\varepsilon _{2}E_{\eta }^{1}\ Curl%
\overrightarrow{B}.\overrightarrow{T}-\kappa E_{s}^{1}-E_{\xi }^{3}\Gamma
_{NB}^{\xi }=0\Rightarrow \kappa =\varepsilon _{2}\frac{E_{\eta }^{1}}{%
E_{s}^{1}}\ Curl \overrightarrow{B}.\overrightarrow{T}+\frac{E_{\xi }^{3}%
}{E_{s}^{1}}\ Div \overrightarrow{B}.  \tag{4.14}
\end{equation}

Thus, for $\overrightarrow{E}$ the following derivative equations can be
written as
\begin{equation}
\frac{\partial \overrightarrow{E}^{s}}{\partial s}=-\varepsilon _{1}\kappa
E_{s}^{1}\overrightarrow{T}-\varepsilon _{2}E_{s}^{3}\tau \overrightarrow{N}%
+\varepsilon _{3}E_{s}^{1}\tau \overrightarrow{B}  \tag{4.15a}
\end{equation}%
\begin{equation}
\frac{\partial \overrightarrow{E}^{\xi }}{\partial \xi }=-\varepsilon
_{1}\left( E_{\xi }^{1}\Gamma _{TN}^{\xi }+E_{\xi }^{3}\Gamma _{TB}^{\xi
}\right) \overrightarrow{T}-\varepsilon _{2}E_{\xi }^{3}\Gamma _{NB}^{\xi }%
\overrightarrow{N}+\varepsilon _{3}E_{\xi }^{1}\Gamma _{NB}^{\xi }%
\overrightarrow{B}  \tag{4.15b}
\end{equation}%
\begin{equation}
\frac{\partial \overrightarrow{E}^{\eta }}{\partial \eta }=-\varepsilon
_{1}\left( E_{\eta }^{1}\Upsilon _{TN}^{\eta }+E_{\eta }^{1}\Upsilon
_{TB}^{\eta }\right) \overrightarrow{T}-\varepsilon _{2}E_{\eta
}^{3}\Upsilon _{NB}^{\eta }\overrightarrow{N}+\varepsilon _{3}E_{\eta
}^{1}\Upsilon _{NB}^{\eta }\overrightarrow{B}.  \tag{4.15c}
\end{equation}

When the particle is affected by the electromagnetic field in the $\xi -$%
direction for the first case, a Lorentz force $\phi _{\xi }$ arises and the
particle moves along a new electromagnetic trajectory according to the frame
in Space form. The electromagnetic vector field $\overrightarrow{M}^{_{\xi }}
$ of the curve $\gamma $ in the $\xi -$direction of the optical fiber for
the first case with respect to the frame satisfies the following condition%
\begin{equation}
\phi _{\xi }(\overrightarrow{E})=\frac{\partial \overrightarrow{E}}{\partial
\xi }=\overrightarrow{M}^{_{\xi }}\times \overrightarrow{E},  \tag{4.16}
\end{equation}%
Lorentz force equation $\phi _{\xi }$ in the $\xi -$direction of the optical
fiber with respect to the frame can be obtain. Hence, by using (4.2) the
derivative equation for $\overrightarrow{E}$ in the $\xi -$direction can be
written as follows 
\begin{equation*}
\frac{\partial \overrightarrow{E}}{\partial \xi }=-\varepsilon
_{1}(\varepsilon _{2}\Gamma _{TN}^{\xi }\overrightarrow{E}\cdot 
\overrightarrow{N}+\varepsilon _{3}\Gamma _{TB}^{\xi }\overrightarrow{E}%
\cdot \overrightarrow{B})\overrightarrow{T}+\varepsilon _{2}(\varepsilon
_{1}\Gamma _{TN}^{\xi }\overrightarrow{E}\cdot \overrightarrow{T}+\ Div%
\overrightarrow{B}\overrightarrow{E}\cdot \overrightarrow{B})\overrightarrow{%
N}
\end{equation*}%
\begin{equation}
+\varepsilon _{3}(\varepsilon _{1}\Gamma _{TB}^{\xi }\overrightarrow{E}\cdot 
\overrightarrow{T}-\ Div \overrightarrow{B}.\overrightarrow{E}\cdot 
\overrightarrow{N})\overrightarrow{B},  \tag{4.17}
\end{equation}

Now, when we consider the components $c_{i}^{\xi },$ $i=1,2,3$ in equation
(4.2) together with $\overrightarrow{E}$, it can be obtained as follows 
\begin{eqnarray*}
\varepsilon _{1}c_{1}^{\xi } &=&\overrightarrow{T}\cdot \frac{\partial 
\overrightarrow{E}}{\partial \xi }=-\overrightarrow{E}.\left( \varepsilon
_{2}\Gamma _{TN}^{\xi }\overrightarrow{N}+\varepsilon _{3}\Gamma _{TB}^{\xi }%
\overrightarrow{B}\right)  \\
\varepsilon _{2}c_{2}^{\xi } &=&\overrightarrow{N}\cdot \frac{\partial 
\overrightarrow{E}}{\partial \xi }=-\overrightarrow{E}.\left( -\varepsilon
_{1}\Gamma _{TN}^{\xi }\overrightarrow{T}-\ Div \overrightarrow{B}\cdot 
\overrightarrow{B}\right)  \\
\varepsilon _{3}c_{3}^{\xi } &=&\overrightarrow{B}\cdot \frac{\partial 
\overrightarrow{E}}{\partial \xi }=-\overrightarrow{E}.\left( -\varepsilon
_{1}\Gamma _{TB}^{\xi }\overrightarrow{T}+\ Div \overrightarrow{B}\cdot 
\overrightarrow{N}\right) .
\end{eqnarray*}

Then, from previous equations and (4.16), we get%
\begin{eqnarray*}
\phi _{\xi }\left( T\right)  &=&\varepsilon _{2}\Gamma _{TN}^{\xi }%
\overrightarrow{N}+\varepsilon _{3}\Gamma _{TB}^{\xi }\overrightarrow{B} \\
\phi _{\xi }\left( N\right)  &=&-\varepsilon _{1}\Gamma _{TN}^{\xi }%
\overrightarrow{T}-\varepsilon _{2}\varepsilon _{3}\ Div \overrightarrow{%
B}\overrightarrow{B} \\
\phi _{\xi }\left( B\right)  &=&-\varepsilon _{1}\Gamma _{TB}^{\xi }%
\overrightarrow{T}+\varepsilon _{2}\varepsilon _{3}\ Div \overrightarrow{%
B}\overrightarrow{N}
\end{eqnarray*}%
also, the Lorentz force equation $\phi _{_{\xi }}$ in the $\xi -$direction
of the optical fiber for the first case with respect to the frame in $%
\mathbf{M}_{q}^{3}$ is written as 
\begin{equation}
\begin{bmatrix}
\phi _{_{\xi }}\left( T\right)  \\ 
\phi _{_{\xi }}\left( N\right)  \\ 
\phi _{_{\xi }}\left( B\right) 
\end{bmatrix}%
=%
\begin{bmatrix}
0 & \varepsilon _{2}\Gamma _{TN}^{\xi } & \varepsilon _{3}\Gamma _{TB}^{\xi }
\\ 
-\varepsilon _{1}\Gamma _{TN}^{\xi } & 0 & -\varepsilon _{2}\varepsilon _{3}%
\ Div \overrightarrow{B} \\ 
-\varepsilon _{1}\Gamma _{TB}^{\xi } & \varepsilon _{2}\varepsilon _{3}\ Div\overrightarrow{B} & 0%
\end{bmatrix}%
\begin{bmatrix}
T \\ 
N \\ 
B%
\end{bmatrix}%
.  \tag{4.18}
\end{equation}

A electromagnetic curve $\gamma $ of the electromagnetic wave in the $\xi -$%
direction along the optical fiber is a magnetic trajectory of a magnetic
field $\overrightarrow{M}^{_{\xi }}$ according to the frame $\left\{
T,N,B\right\} $ in $\mathbf{M}_{q}^{3}$ and this magnetic field $%
\overrightarrow{M}^{_{\xi }}$ is obtained as%
\begin{equation}
\overrightarrow{M}^{_{\xi }}=m_{1}^{\xi }\overrightarrow{T}+m_{2}^{\xi }%
\overrightarrow{N}+m_{3}^{\xi }\overrightarrow{B},  \tag{4.19}
\end{equation}%
where $m_{i}^{\xi },$ $i=1,2,3$ are smooth functions. The following system
of equations is obtained from equation (4.16), (4.19) and (3.21)%
\begin{equation}
\overrightarrow{M}^{_{\xi }}\times \overrightarrow{T}=\phi _{_{\xi }}(%
\overrightarrow{T})=\frac{\partial \overrightarrow{T}}{\partial \xi }%
=-\varepsilon _{3}m_{2}^{\xi }\overrightarrow{B}+\varepsilon _{2}m_{3}^{\xi }%
\overrightarrow{N}=\varepsilon _{2}\Gamma _{TN}^{\xi }\overrightarrow{N}%
+\varepsilon _{3}\Gamma _{TB}^{\xi }\overrightarrow{B}  \tag{4.20a}
\end{equation}%
\begin{equation}
\overrightarrow{M}^{_{\xi }}\times \overrightarrow{N}=\phi _{_{\xi }}(%
\overrightarrow{N})=\frac{\partial \overrightarrow{N}}{\partial \xi }%
=\varepsilon _{3}m_{1}^{\xi }\overrightarrow{B}-\varepsilon _{1}m_{3}^{\xi }%
\overrightarrow{T}=-\varepsilon _{1}\Gamma _{TN}^{\xi }\overrightarrow{T}%
-\varepsilon _{2}\varepsilon _{3}\ Div \overrightarrow{B}\overrightarrow{%
B}  \tag{4.20b}
\end{equation}%
\begin{equation}
\overrightarrow{M}^{_{\xi }}\times \overrightarrow{B}=\phi _{_{\xi }}(%
\overrightarrow{B})=\frac{\partial \overrightarrow{B}}{\partial \xi }%
=-\varepsilon _{2}m_{1}^{\xi }\overrightarrow{N}+\varepsilon _{1}m_{2}^{\xi }%
\overrightarrow{T}=-\varepsilon _{1}\Gamma _{TB}^{\xi }\overrightarrow{T}%
+\varepsilon _{2}\varepsilon _{3}\ Div \overrightarrow{B}\overrightarrow{%
N}.  \tag{4.20c}
\end{equation}

In the above equation system, the coefficients are found as follows, taking
into account the algebraic equations. 
\begin{eqnarray*}
-m_{2}^{\xi } &=&\Gamma _{TB}^{\xi },\text{ }m_{3}^{\xi }=\Gamma _{TN}^{\xi }%
\text{; }m_{1}^{\xi }=-\varepsilon _{2}\ Div \overrightarrow{B}, \\
\text{ }-m_{3}^{\xi } &=&-\Gamma _{TN}^{\xi }\text{; }-m_{1}^{\xi
}=\varepsilon _{3}\ Div \overrightarrow{B},\text{ }m_{2}^{\xi }=-\Gamma
_{TB}^{\xi }
\end{eqnarray*}%
and we get%
\begin{equation}
\overrightarrow{M}^{_{\xi }}=-\varepsilon _{2}\ Div \overrightarrow{B}%
\overrightarrow{T}-\Gamma _{TB}^{\xi }\overrightarrow{N}+\Gamma _{TN}^{\xi }%
\text{ }\overrightarrow{B}.  \tag{4.21}
\end{equation}

If the derivative with respect to $s$ is taken in (4.21) and the inner
product with $\overrightarrow{T}$ is made, the following equation is obtained%
\begin{equation}
\overrightarrow{T}\cdot \frac{\partial \overrightarrow{M}^{_{\xi }}}{%
\partial s}=\overrightarrow{T}\cdot \left( 
\begin{array}{c}
-\varepsilon _{2}\frac{\partial \ Div \overrightarrow{B}}{\partial s}%
\overrightarrow{T}-\varepsilon _{2}\ Div \overrightarrow{B}\frac{%
\partial \overrightarrow{T}}{\partial s}-\frac{\partial \Gamma _{TB}^{\xi }}{%
\partial s}\overrightarrow{N} \\ 
-\Gamma _{TB}^{\xi }\frac{\partial \overrightarrow{N}}{\partial s}+\frac{%
\partial \Gamma _{TN}^{\xi }}{\partial s}\overrightarrow{B}+\Gamma
_{TN}^{\xi }\text{ }\frac{\partial \overrightarrow{B}}{\partial s}%
\end{array}%
\right) .  \tag{4.22}
\end{equation}

Finally, if $\frac{\partial \overrightarrow{T}}{\partial s},\frac{\partial 
\overrightarrow{N}}{\partial s},\frac{\partial \overrightarrow{B}}{\partial s%
}$ are written in the last equation, we get%
\begin{equation}
\overrightarrow{T}\cdot \frac{\partial \overrightarrow{M}^{_{\xi }}}{%
\partial s}=-\varepsilon _{2}\varepsilon _{1}\frac{\partial \ Div %
\overrightarrow{B}}{\partial s}-\kappa \Gamma _{TB}^{\xi }.  \tag{4.23}
\end{equation}

Similarly, firstly using the equations (3.8)( or (3.21)) and (3.9)(or
(3.22)) respectively, the following equations are obtained for $%
\overrightarrow{N}$ and $\overrightarrow{B}$%
\begin{equation}
\overrightarrow{N}\cdot \frac{\partial \overrightarrow{M}^{_{\xi }}}{%
\partial \xi }=-\varepsilon _{3}\ Div \overrightarrow{B}\Gamma
_{TN}^{\xi }-\varepsilon _{2}\frac{\partial \Gamma _{TB}^{\xi }}{\partial
\xi }-\Gamma _{NB}^{\xi }\Gamma _{TN}^{\xi }\text{ }  \tag{4.24}
\end{equation}%
\begin{equation}
\overrightarrow{B}\cdot \frac{\partial \overrightarrow{M}^{_{\xi }}}{%
\partial \eta }=-\varepsilon _{2}\ Div \overrightarrow{B}\Upsilon
_{TB}^{\eta }-\Gamma _{TB}^{\xi }\Upsilon _{NB}^{\eta }+\varepsilon _{3}%
\frac{\partial \Gamma _{TB}^{\xi }}{\partial \eta }.  \tag{4.25}
\end{equation}

Considering the Maxwell equations and using the equations (4,23), (4.24),
(4,25), the following expression is obtained%
\begin{equation*}
\nabla \overrightarrow{M}^{_{\xi }}=-\varepsilon _{2}\varepsilon _{1}\frac{%
\partial \ Div \overrightarrow{B}}{\partial s}-\kappa \Gamma _{TB}^{\xi
}-\varepsilon _{3}\ Div \overrightarrow{B}\Gamma _{TN}^{\xi }
\end{equation*}%
\begin{equation}
-\varepsilon _{2}\frac{\partial \Gamma _{TB}^{\xi }}{\partial \xi }-\Gamma
_{NB}^{\xi }\Gamma _{TN}^{\xi }-\varepsilon _{2}\ Div \overrightarrow{B}%
\Upsilon _{TB}^{\eta }-\Gamma _{TB}^{\xi }\Upsilon _{NB}^{\eta }+\varepsilon
_{3}\frac{\partial \Gamma _{TB}^{\xi }}{\partial \eta }.  \tag{4.26}
\end{equation}

Since equality is equal to zero in the Maxwell equations, the following
equation can be written
\begin{equation}
\kappa =\frac{-1}{\Gamma _{TB}^{\xi }}\left( 
\begin{array}{c}
\varepsilon _{2}\varepsilon _{1}\frac{\partial \ Div \overrightarrow{B}}{%
\partial s}+\varepsilon _{3}\ Div \overrightarrow{B}\Gamma _{TN}^{\xi
}+\varepsilon _{2}\frac{\partial \Gamma _{TB}^{\xi }}{\partial \xi }+\Gamma
_{NB}^{\xi }\Gamma _{TN}^{\xi } \\ 
+\varepsilon _{2}\ Div \overrightarrow{B}\Upsilon _{TB}^{\eta }+\Gamma
_{TB}^{\xi }\Upsilon _{NB}^{\eta }-\varepsilon _{3}\frac{\partial \Gamma
_{TB}^{\xi }}{\partial \eta }%
\end{array}%
\right) .  \tag{4.27}
\end{equation}

Moreover, if we consider that the electric field is right angle to the
tangential direction and by taking the derivatives of the vector field
defined in (4.21) with respect to $s$, $\xi $, $\eta $, respectively, we get
\begin{equation}
\nabla \times \overrightarrow{M}^{\xi }=\overrightarrow{T}\times \frac{%
\partial \overrightarrow{M}^{_{\xi }}}{\partial s}+\overrightarrow{N}\times 
\frac{\partial \overrightarrow{M}^{_{\xi }}}{\partial \xi }+\overrightarrow{B%
}\times \frac{\partial \overrightarrow{M}^{_{\xi }}}{\partial \eta }. 
\tag{4.28}
\end{equation}

When the calculation is made for the three values {}{}in the previous
equation, the following equations are obtained%
\begin{eqnarray*}
\overrightarrow{T}\times \frac{\partial \overrightarrow{M}^{_{\xi }}}{%
\partial s} &=&-\varepsilon _{2}\varepsilon _{1}c\ Div \overrightarrow{B}%
\overrightarrow{T}\times \overrightarrow{\gamma }-\varepsilon
_{2}(\varepsilon _{3}\tau \Gamma _{TB}^{\xi }+\frac{\partial \Gamma
_{TN}^{\xi }}{\partial s})\overrightarrow{N} \\
&&-\varepsilon _{3}(\ Div \overrightarrow{B}\kappa +\frac{\partial
\Gamma _{TB}^{\xi }}{\partial s}+\varepsilon _{2}\tau \Gamma _{TN}^{\xi })%
\overrightarrow{B} \\
\overrightarrow{N}\times \frac{\partial \overrightarrow{M}^{_{\xi }}}{%
\partial \xi } &=&\varepsilon _{1}(\frac{\partial \Gamma _{TN}^{\xi }}{%
\partial \xi }-\varepsilon _{2}\varepsilon _{3}\Gamma _{TB}^{\xi }\ Div %
\overrightarrow{B}-\varepsilon _{3}\Gamma _{TB}^{\xi }\Gamma _{NB}^{\xi })%
\overrightarrow{T}+\varepsilon _{2}\varepsilon _{3}\frac{\partial \ Div %
\overrightarrow{B}}{\partial \xi }\overrightarrow{B} \\
\overrightarrow{B}\times \frac{\partial \overrightarrow{M}^{_{\xi }}}{%
\partial \eta } &=&\varepsilon _{1}(\ Div \overrightarrow{B}\Upsilon
_{TN}^{\eta }+\frac{\partial \Gamma _{TB}^{\xi }}{\partial \eta }%
+\varepsilon _{2}\Gamma _{TN}^{\xi }\Upsilon _{NB}^{\eta })\overrightarrow{T}
\\
&&+(-\frac{\partial \ Div \overrightarrow{B}}{\partial \eta }%
+\varepsilon _{1}\varepsilon _{2}\Gamma _{TB}^{\xi }\Upsilon _{TN}^{\eta
} -\varepsilon _{1}\varepsilon _{2}\Gamma _{TN}^{\xi }\Upsilon _{TB}^{\eta })%
\overrightarrow{N}
\end{eqnarray*}%
and by using these equations, we write 
\begin{eqnarray*}
\nabla \times \overrightarrow{M}^{_{\xi }} &=&-\varepsilon _{2}\varepsilon
_{1}c\ Div B\overrightarrow{T}\times \overrightarrow{\gamma }%
+\varepsilon _{1}\left( 
\begin{array}{c}
\frac{\partial \Gamma _{TN}^{\xi }}{\partial \xi }-\varepsilon
_{2}\varepsilon _{3}\Gamma _{TB}^{\xi }\ Div \overrightarrow{B}%
+\varepsilon _{2}\Gamma _{TN}^{\xi }\Upsilon _{NB}^{\eta } \\ 
-\varepsilon _{3}\Gamma _{TB}^{\xi }\Gamma _{NB}^{\xi }\Upsilon _{TN}^{\eta }%
\ Div \overrightarrow{B}+\frac{\partial \Gamma _{TB}^{\xi }}{\partial
\eta }%
\end{array}%
\right) \overrightarrow{T} \\
&&+\left( -\varepsilon _{2}\varepsilon _{3}\tau \Gamma _{TB}^{\xi
}-\varepsilon _{2}\frac{\partial \Gamma _{TN}^{\xi }}{\partial s}-\frac{%
\partial \ Div \overrightarrow{B}}{\partial \eta }+\varepsilon
_{1}\varepsilon _{2}\Gamma _{TB}^{\xi }\Upsilon _{TN}^{\eta }-\varepsilon
_{1}\varepsilon _{2}\Gamma _{TN}^{\xi }\Upsilon _{TB}^{\eta }\right) 
\overrightarrow{N}
\end{eqnarray*}%
\begin{equation}
+\varepsilon _{3}\left( -\kappa \ Div \overrightarrow{B}-\frac{\partial
\Gamma _{TB}^{\xi }}{\partial s}-\varepsilon _{2}\tau \Gamma _{TN}^{\xi
}+\varepsilon _{2}\varepsilon _{3}\frac{\partial \ Div \overrightarrow{B}%
}{\partial \xi }\right) \overrightarrow{B}.  \tag{4.29}
\end{equation}

As the second situation, Lorentz force equation $\phi _{\eta }$ in the $\eta
-$direction of the optical fiber with respect to the frame can be obtain. By
performing similar algebraic calculations 
\begin{equation}
\frac{\partial \overrightarrow{E}}{\partial \eta }=b_{1}^{\eta }%
\overrightarrow{T}+b_{2}^{\eta }\overrightarrow{N}+b_{3}^{\eta }%
\overrightarrow{B},  \tag{4.30}
\end{equation}%
where $b_{i}^{\eta },$ $i=1,2,3$ are smooth functions. From (3.9) we obtain 
\begin{eqnarray*}
\varepsilon _{1}b_{1}^{\eta } &=&-\overrightarrow{E}_{_{s\xi \eta }}\cdot
\left( \varepsilon _{2}\Upsilon _{TN}^{\eta }\overrightarrow{N}+\varepsilon
_{3}\Upsilon _{TB}^{\eta }\overrightarrow{B}\right) ;\varepsilon
_{2}b_{2}^{\eta }=-\overrightarrow{E}_{_{s\xi \eta }}\cdot \left(
-\varepsilon _{1}\Upsilon _{TN}^{\eta }\overrightarrow{T}+\varepsilon
_{3}\Upsilon _{NB}^{\eta }\overrightarrow{B}\right) ; \\
\varepsilon _{3}b_{3}^{\eta } &=&\overrightarrow{E}_{_{s\xi \eta }}\cdot
\left( \varepsilon _{1}\Upsilon _{TB}^{\eta }\overrightarrow{T}+\varepsilon
_{2}\Upsilon _{NB}^{\eta }\overrightarrow{N}\right) .
\end{eqnarray*}

Considering the last equations in (4.30), we have%
\begin{eqnarray*}
\frac{\partial \overrightarrow{E}}{\partial \eta } &=&-\varepsilon
_{1}(\varepsilon _{2}\Upsilon _{TN}^{\eta }\overrightarrow{E}\cdot 
\overrightarrow{N}+\varepsilon _{3}\Upsilon _{TB}^{\eta }\overrightarrow{E}%
\cdot \overrightarrow{B})\overrightarrow{T} \\
&&+\varepsilon _{2}(\varepsilon _{1}\Upsilon _{TN}^{\eta }\overrightarrow{E}%
\cdot \overrightarrow{T}-\varepsilon _{2}\Upsilon _{NB}^{\eta }%
\overrightarrow{E}\cdot \overrightarrow{B})\overrightarrow{N}+\varepsilon
_{3}(\varepsilon _{1}\Upsilon _{TB}^{\eta }\overrightarrow{E}\cdot 
\overrightarrow{T}+\varepsilon _{2}\Upsilon _{NB}^{\eta }\overrightarrow{E}%
\cdot \overrightarrow{N})\overrightarrow{B}.
\end{eqnarray*}

Then,\ by using previous equations and (3.9) we can obtain the equation
given as 
\begin{equation}
\phi _{\eta }(\overrightarrow{E})=\frac{\partial \overrightarrow{E}}{%
\partial \eta }=\overrightarrow{M}^{_{^{\eta }}}\times \overrightarrow{E}, 
\tag{4.31}
\end{equation}%
we say that the electromagnetic vector field $\overrightarrow{M}^{_{^{\eta
}}}$ of the curve $\gamma $ in the $\eta -$direction of the optical fiber
for the second case with respect to the frame satisfies the this equation.
Also, by using (4.31) we can calculate%
\begin{eqnarray*}
\phi _{\eta }\left( T\right)  &=&\frac{\partial \overrightarrow{T}}{\partial
\eta }=\varepsilon _{2}\Upsilon _{TN}^{\eta }\overrightarrow{N}+\varepsilon
_{3}\Upsilon _{TB}^{\eta }\overrightarrow{B} \\
\phi _{\eta }\left( N\right)  &=&\frac{\partial \overrightarrow{N}}{\partial
\eta }=-\varepsilon _{1}\Upsilon _{TN}^{\eta }\overrightarrow{T}+\varepsilon
_{3}\Upsilon _{NB}^{\eta }\overrightarrow{B} \\
\phi _{\eta }\left( B\right)  &=&\frac{\partial \overrightarrow{B}}{\partial
\eta }=-\varepsilon _{1}\Upsilon _{TB}^{\eta }\overrightarrow{T}-\varepsilon
_{2}\Upsilon _{NB}^{\eta }\overrightarrow{N}.
\end{eqnarray*}

Therefore, the Lorentz force equation $\phi _{\eta }$ in the $\eta -$%
direction of the optical fiber with respect to the frame in $\mathbf{M}%
_{q}^{3}$ is written as 
\begin{equation}
\begin{bmatrix}
\phi _{\eta }\left( T\right)  \\ 
\phi _{\eta }\left( N\right)  \\ 
\phi _{\eta }\left( B\right) 
\end{bmatrix}%
=%
\begin{bmatrix}
0 & \varepsilon _{2}\Upsilon _{TN}^{\eta } & \varepsilon _{3}\Upsilon
_{TB}^{\eta } \\ 
-\varepsilon _{1}\Upsilon _{TN}^{\eta } & 0 & \varepsilon _{3}\Upsilon
_{NB}^{\eta } \\ 
-\varepsilon _{1}\Upsilon _{TB}^{\eta } & -\varepsilon _{2}\Upsilon
_{NB}^{\eta } & 0%
\end{bmatrix}%
\begin{bmatrix}
T \\ 
N \\ 
B%
\end{bmatrix}%
.  \tag{4.32}
\end{equation}

A electromagnetic curve $\gamma $ of the electromagnetic wave in the $\eta -$%
direction along the optical fiber is a magnetic trajectory of a magnetic
field $\overrightarrow{M}^{_{^{\eta }}}$ according to the frame $\left\{
T,N,B\right\} $ in $\mathbf{M}_{q}^{3}$ and this magnetic field $%
\overrightarrow{M}^{_{^{\eta }}}$ is obtained as 
\begin{equation}
\overrightarrow{M}^{_{^{\eta }}}=m_{1}^{\eta }\overrightarrow{T}+m_{2}^{\eta
}\overrightarrow{N}+m_{3}^{\eta }\overrightarrow{B},  \tag{4.33}
\end{equation}%
where $m_{i}^{\eta },$ $i=1,2,3$ are smooth functions. Hence, from (4.31)
and (4.33) we get%
\begin{eqnarray*}
\overrightarrow{M}^{_{^{\eta }}}\times \overrightarrow{T} &=&\phi _{\eta
}(T)=-\varepsilon _{3}m_{2}^{\eta }\overrightarrow{B}+\varepsilon
_{2}m_{3}^{\eta }\overrightarrow{N}=\varepsilon _{2}\Upsilon _{TN}^{\eta }%
\overrightarrow{N}+\varepsilon _{3}\Upsilon _{TB}^{\eta }\overrightarrow{B}
\\
\overrightarrow{M}^{_{^{\eta }}}\times \overrightarrow{N} &=&\phi _{\eta
}(N)=\varepsilon _{3}m_{1}^{\eta }\overrightarrow{B}-\varepsilon
_{1}m_{3}^{\eta }\overrightarrow{T}=-\varepsilon _{1}\Upsilon _{TN}^{\eta }%
\overrightarrow{T}+\varepsilon _{3}\Upsilon _{NB}^{\eta }\overrightarrow{B}
\\
\overrightarrow{M}^{_{^{\eta }}}\times \overrightarrow{B} &=&\phi _{\eta
}(B)=-\varepsilon _{2}m_{1}^{\eta }\overrightarrow{N}+\varepsilon
_{1}m_{2}^{\eta }\overrightarrow{T}=-\varepsilon _{1}\Upsilon _{TB}^{\eta }%
\overrightarrow{T}-\varepsilon _{2}\Upsilon _{NB}^{\eta }\overrightarrow{N}
\end{eqnarray*}%
and by taking into account the algebraic equations the coefficients are
found as follows 
\begin{equation*}
-m_{2}^{\eta }=\Upsilon _{TB}^{\eta },m_{3}^{\eta }=\Upsilon _{TN}^{\eta }%
\text{;}m_{1}^{\eta }=\Upsilon _{NB}^{\eta },m_{3}^{\eta }=\Upsilon
_{TN}^{\eta }\text{;}m_{1}^{\eta }=-\Upsilon _{NB}^{\eta },m_{2}^{\eta
}=-\Upsilon _{TB}^{\eta }.
\end{equation*}

Hence, from (3.22) we get%
\begin{equation}
\overrightarrow{M}^{_{^{\eta }}}=\varepsilon _{2}(\ Curl \overrightarrow{%
B}\cdot \overrightarrow{T})\overrightarrow{T}+\varepsilon _{1}(\varepsilon
_{2}\tau +\ Curl \overrightarrow{B}\cdot \overrightarrow{N})%
\overrightarrow{N}-\varepsilon _{1}(\varepsilon _{3}\tau +\ Curl %
\overrightarrow{N}\cdot \overrightarrow{N})\overrightarrow{B}.  \tag{4.34}
\end{equation}

If the derivative is taken with respect to $s$ in (4.34) and the derivative
equations given in (2.2) are taken into account, for $\overrightarrow{T}$
the following equation is obtained
\begin{equation}
\overrightarrow{T}\cdot \frac{\partial \overrightarrow{M}^{_{^{\eta }}}}{%
\partial s}=\varepsilon _{1}\varepsilon _{2}\frac{\partial }{\partial s}%
\left( \ Curl \overrightarrow{B}\cdot \overrightarrow{T}\right) +\kappa
\left( \varepsilon _{2}\tau +\ Curl \overrightarrow{B}\cdot 
\overrightarrow{N}\right) .  \tag{4.35}
\end{equation}

Similarly, if the derivatives are taken with respect to $\xi $ and $\eta $
in (4.34), respectively, and the derivative equations given in (3.21) and
(3.22) are taken into account, for $\overrightarrow{N}$ and $\overrightarrow{%
B}$ the following equations are obtained%
\begin{equation}
\overrightarrow{N}\cdot \frac{\partial \overrightarrow{M}^{_{^{\eta }}}}{%
\partial \xi }=\varepsilon _{2}\varepsilon _{1}\frac{\partial (\varepsilon
_{2}\tau +\ Curl \overrightarrow{B}.\overrightarrow{N})}{\partial \xi }%
\text{ }+\varepsilon _{2}\ Curl \overrightarrow{B}.\overrightarrow{T}%
\Gamma _{TN}^{\xi }+\varepsilon _{1}(\varepsilon _{3}\tau +\ Curl 
\overrightarrow{N}.\overrightarrow{N})\Gamma _{NB}^{\xi }  \tag{4.36}
\end{equation}%
\begin{equation}
\overrightarrow{B}\cdot \frac{\partial \overrightarrow{M}^{_{^{\eta }}}}{%
\partial \eta }=\varepsilon _{2}\ Curl \overrightarrow{B}.%
\overrightarrow{T}\Upsilon _{TB}^{\eta }+\varepsilon _{1}(\varepsilon
_{2}\tau +\ Curl \overrightarrow{B}.\overrightarrow{N})\Upsilon
_{NB}^{\eta }-\varepsilon _{3}\varepsilon _{1}\frac{\partial (\varepsilon
_{3}\tau +\ Curl\overrightarrow{N}.\overrightarrow{N})}{\partial \eta }.
\tag{4.37}
\end{equation}

Finally, from Maxwell's equations and equations (4.35), (4.36) and (4.37),
the following equality is obtained%
\begin{equation*}
\nabla \overrightarrow{M}^{_{^{\eta }}}=\varepsilon _{1}\varepsilon _{2}%
\frac{\partial \ Curl \overrightarrow{B}.\overrightarrow{T}}{\partial s}%
+\kappa (\varepsilon _{2}\tau +\ Curl \overrightarrow{B}.\overrightarrow{%
N})+\varepsilon _{2}\varepsilon _{1}\frac{\partial (\varepsilon _{2}\tau +%
\ Curl \overrightarrow{B}.\overrightarrow{N})}{\partial \xi }\text{ }
\end{equation*}%
\begin{equation*}
+\varepsilon _{2}\ Curl \overrightarrow{B}.\overrightarrow{T}\Gamma
_{TN}^{\xi }+\varepsilon _{1}(\varepsilon _{3}\tau +\ Curl %
\overrightarrow{N}.\overrightarrow{N})\Gamma _{NB}^{\xi }+\varepsilon _{2}%
\ Curl \overrightarrow{B}.\overrightarrow{T}\Upsilon _{TB}^{\eta }
\end{equation*}%
\begin{equation}
+\varepsilon _{1}(\varepsilon _{2}\tau +\ Curl \overrightarrow{B}.%
\overrightarrow{N})\Upsilon _{NB}^{\eta }-\varepsilon _{3}\varepsilon _{1}%
\frac{\partial (\varepsilon _{3}\tau +\ Curl \overrightarrow{N}.%
\overrightarrow{N})}{\partial \eta }  \tag{4.38}
\end{equation}%
from previous equation, we have
\begin{equation}
\kappa =\frac{-1}{\varepsilon _{2}\tau +\ Curl \overrightarrow{B}.%
\overrightarrow{N}}\left( 
\begin{array}{c}
\varepsilon _{1}\varepsilon _{2}\frac{\partial (\ Curl \overrightarrow{B}%
.\overrightarrow{T})}{\partial s}+\varepsilon _{2}\varepsilon _{1}\frac{%
\partial (\varepsilon _{2}\tau +\ Curl \overrightarrow{B}.%
\overrightarrow{N})}{\partial \xi }\text{ }+\varepsilon _{2}\ Curl %
\overrightarrow{B}.\overrightarrow{T}\Gamma _{TN}^{\xi } \\ 
-\varepsilon _{3}\varepsilon _{1}\frac{\partial (\varepsilon _{3}\tau +\ Curl \overrightarrow{N}.\overrightarrow{N})}{\partial \eta }+\varepsilon
_{1}(\varepsilon _{3}\tau +\ Curl \overrightarrow{N}.\overrightarrow{N}%
)\Gamma _{NB}^{\xi } \\ 
+\varepsilon _{2}\ Curl \overrightarrow{B}.\overrightarrow{T}\Upsilon
_{TB}^{\eta }+\varepsilon _{1}(\varepsilon _{2}\tau +\ Curl %
\overrightarrow{B}.\overrightarrow{N})\Upsilon _{NB}^{\eta }%
\end{array}%
\right) .  \tag{4.39}
\end{equation}

Similarly, if partial derivatives are taken with respect to $s$, $\xi $, $%
\eta $ in the expression given in (4.34) and used in following equation
given as
\begin{equation}
\nabla \times \overrightarrow{M}^{_{^{\eta }}}=\overrightarrow{T}\times 
\frac{\partial \overrightarrow{M}^{_{^{\eta }}}}{\partial s}+\overrightarrow{%
N}\times \frac{\partial \overrightarrow{M}^{_{^{\eta }}}}{\partial \xi }+%
\overrightarrow{B}\times \frac{\partial \overrightarrow{M}^{_{^{\eta }}}}{%
\partial \eta }  \tag{4.40}
\end{equation}%
and from (4.40), we get 
\begin{equation*}
\nabla \times \overrightarrow{M}^{_{^{\eta }}}=-\varepsilon _{2}\varepsilon
_{1}c\ Curl \overrightarrow{B}.\overrightarrow{T}\overrightarrow{T}%
\times \overrightarrow{\alpha }+\Theta _{1}^{\eta }\overrightarrow{T}+\Theta
_{2}^{\eta }\overrightarrow{N}+\Theta _{3}^{\eta }\overrightarrow{B},
\end{equation*}%
where 
\begin{eqnarray*}
\Theta _{1}^{\eta } &=&-\varepsilon _{2}\varepsilon _{1}c\ Curl %
\overrightarrow{B}.\overrightarrow{T}\left( \overrightarrow{T}\times 
\overrightarrow{\gamma }\right) +\varepsilon _{1}\varepsilon _{2}\left( 
\begin{array}{c}
\frac{\partial (\varepsilon _{3}\tau +\ Curl\overrightarrow{N}.%
\overrightarrow{N})}{\partial s} \\ 
-\varepsilon _{3}\tau (\varepsilon _{2}\tau +\ Curl \overrightarrow{B}.%
\overrightarrow{N})%
\end{array}%
\right) \overrightarrow{N} \\
&&+\varepsilon _{1}\varepsilon _{3}\left( \varepsilon _{3}\kappa \ Curl%
\overrightarrow{B}.\overrightarrow{T}+\frac{\partial (\varepsilon _{2}\tau +%
\ Curl \overrightarrow{B}.\overrightarrow{N})}{\partial s}+\varepsilon
_{2}(\varepsilon _{3}\tau +\ Curl \overrightarrow{N}.\overrightarrow{N}%
)\right) \overrightarrow{B} \\
\Theta _{2}^{\eta } &=&\left( 
\begin{array}{c}
\varepsilon _{1}\varepsilon _{2}\varepsilon _{3}\ Curl \overrightarrow{B}%
.\overrightarrow{T}\Gamma _{TB}^{\xi }+\varepsilon _{3}(\varepsilon _{2}\tau
+\ Curl \overrightarrow{B}.\overrightarrow{N})\Gamma _{NB}^{\xi } \\ 
-\frac{\partial (\varepsilon _{3}\tau +\ Curl \overrightarrow{N}.%
\overrightarrow{N})}{\partial \xi }%
\end{array}%
\right) \overrightarrow{T} \\
&&-\varepsilon _{3}\left( \varepsilon _{2}\frac{\partial \ Curl%
\overrightarrow{B}.\overrightarrow{T}}{\partial \xi }+(\varepsilon _{2}\tau +%
\ Curl \overrightarrow{B}.\overrightarrow{N})\Gamma _{TN}^{\xi
}+(\varepsilon _{3}\tau +\ Curl\overrightarrow{N}.\overrightarrow{N}%
)\Gamma _{TB}^{\xi }\right) \overrightarrow{B} \\
\Theta _{3}^{\eta } &=&-\left( \varepsilon _{1}\ Curl \overrightarrow{B}.%
\overrightarrow{T}\Upsilon _{TN}^{\eta }+\frac{\partial (\varepsilon
_{3}\tau +\ Curl\overrightarrow{B}.\overrightarrow{N})}{\partial \eta }%
+\varepsilon _{2}(\varepsilon _{3}\tau +\ Curl \overrightarrow{N}.%
\overrightarrow{N})\Gamma _{NB}^{\xi }\right) \overrightarrow{T} \\
&&+\left( 
\begin{array}{c}
\varepsilon _{2}(\varepsilon _{3}\tau +\ Curl \overrightarrow{N}.%
\overrightarrow{N})\Upsilon _{TB}^{\eta }+\varepsilon _{2}\frac{\partial 
\ Curl \overrightarrow{B}.\overrightarrow{T}}{\partial \eta } \\ 
-\varepsilon _{2}(\varepsilon _{2}\tau +\ Curl \overrightarrow{B}.%
\overrightarrow{N})\Upsilon _{TN}^{\eta }%
\end{array}%
\right) \overrightarrow{N}.
\end{eqnarray*}

\section{The energy of the vector fields on a particle in $M_{q}^{3}\left(
c\right) $}

In this section, the bending energy formulas for tangent vector of $s$%
-lines( $\xi $-lines, $\eta $-lines respectively) of elastic curve written
by extended Serret-Frenet relations along the curve $\gamma $ are
investigated in $M_{q}^{3}\left( c\right) $.

\subsection{The energy of unit tangent vector of $s-$lines on a moving
particle in $M_{q}^{3}$}

In the subsection, we calculate the energy of the unit tangent vector of $s$
-lines of the curve in $M_{q}^{3}\left( c\right) $ and we also investigate
the bending energy formula for an elastic curve given by extended
Serret-Frenet relations along the curve $\gamma (s,\xi ,\eta )$ in $%
M_{q}^{3}\left( c\right) .$

Let $P$ be a moving particle in $M_{q}^{3}\left( c\right) $ such that it
corresponds to a curve $\gamma (s,\xi ,\eta )$ with parameter $s$, which $s$
is the distance along the $s$-lines of the curve in $s-$direction and
tangent vector of $s$-lines is defined by $\frac{\partial \overrightarrow{T}%
}{\partial s}.$ Hence, by using Sasaki metric and the equations (2.3),
(2.4), (2.5), the energy on the particle in vector field $\frac{\partial 
\overrightarrow{T}}{\partial s}$ can be written as 
\begin{equation*}
energy_{T_{s}}=\frac{1}{2}\int \rho _{s}(dT(T),dT(T))ds
\end{equation*}%
and%
\begin{equation*}
\rho _{s}(dT(T),dT(T))=\rho _{s}\left( T,T\right) +\rho _{s}\left( \nabla
_{T}T,\nabla _{T}T\right) ,
\end{equation*}%
since $\nabla _{T}T=-c\varepsilon _{1}\overrightarrow{\gamma }+\varepsilon
_{2}\kappa \overrightarrow{N},$ we obtain 
\begin{equation}
energy_{T_{s}}=\frac{1}{2}\int \left( \varepsilon _{1}+c^{2}\left\Vert 
\overrightarrow{\gamma }\right\Vert ^{2}+\varepsilon _{2}\kappa ^{2}\right)
ds.  \tag{5.1}
\end{equation}

Also, the energy on the particle in vector field $\frac{\partial N}{\partial
s}$ is written as 
\begin{equation*}
energy_{N_{s}}=\frac{1}{2}\int \rho _{s}(dN(N),dN(N))ds,
\end{equation*}%
since $\nabla _{N}N=-\varepsilon _{1}\kappa \overrightarrow{T}+\varepsilon
_{3}\tau \overrightarrow{B}$, the energy of the vector field $\frac{\partial
N}{\partial s}$is obtain as  
\begin{equation}
energy_{N_{s}}=\frac{1}{2}\int \left( \varepsilon _{2}+\varepsilon
_{1}\kappa ^{2}+\varepsilon _{3}\tau ^{2}\right) ds.  \tag{5.2}
\end{equation}

Similarly, from $\nabla _{B}B=-\varepsilon _{2}\tau \overrightarrow{N}$, the
energy of the vector field $\frac{\partial B}{\partial s}$ is written as, we
get 
\begin{equation}
energy_{B_{s}}=\frac{1}{2}\int \left( \varepsilon _{3}+\varepsilon _{2}\tau
^{2}\right) ds.  \tag{5.3}
\end{equation}

\subsection{The energy of unit tangent vector of $\protect\xi -$lines on a
moving particle in $M_{q}^{3}\left( c\right) $}

In the subsection, we calculate the energy of the unit tangent vector of $%
\xi -$lines of the curve in $M_{q}^{3}\left( c\right) $ and we also
investigate the bending energy formula for an elastic curve given by
extended Serret-Frenet relations along the curve $\gamma (s,\xi ,\eta )$ in $%
M_{q}^{3}\left( c\right) $, which $\xi $ is the distance along the $\xi -$%
lines of the curve in $\xi -$direction and the tangent vector of $\xi -$%
lines is expressed by $\frac{\partial \gamma ^{\prime }}{\partial \xi }.$
Hence, the energy on the particle in vector field $\frac{\partial \gamma
^{\prime }}{\partial \xi }$ can be written as 
\begin{equation*}
energy_{T_{\xi }}=\frac{1}{2}\int \rho _{\xi }(dT(T),dT(T))d\xi ,
\end{equation*}%
from (2.3), (2.4), (2.5), we get%
\begin{equation*}
\rho _{\xi }(dT(T),dT(T))=\rho _{\xi }\left( T,T\right) +\rho _{\xi }\left(
\nabla _{T}T,\nabla _{T}T\right) 
\end{equation*}%
by using the extended Serret-Frenet relations according to parameter $\xi ,$
since $\frac{\partial T}{\partial \xi }=-\varepsilon _{1}\varepsilon _{3}(%
\ Curl \overrightarrow{N}\cdot \overrightarrow{B})\overrightarrow{N}%
+\varepsilon _{1}\varepsilon _{3}\Psi _{B}\overrightarrow{\alpha },$ we get 
\begin{equation}
energy_{T_{\xi }}=\frac{1}{2}\int \left( \varepsilon _{1}+\varepsilon
_{2}\left( \ Curl \overrightarrow{N}.\overrightarrow{B}\right)
^{2}+\varepsilon _{3}\left( \ Curl \overrightarrow{B}.\overrightarrow{B}%
\right) ^{2}\right) d\xi .  \tag{5.4}
\end{equation}

Also, the energy on the particle in vector field $\frac{\partial N}{\partial
\xi }$ is written as 
\begin{equation*}
energy_{N_{\xi }}=\frac{1}{2}\int \rho _{\xi }(dN(N),dN(N))d\xi 
\end{equation*}%
and since $\nabla _{N}N=-(\ Curl \overrightarrow{N}.\overrightarrow{B})%
\overrightarrow{T}+\varepsilon _{3}(-\ Div \overrightarrow{B})%
\overrightarrow{B}$ we can write as%
\begin{equation}
energy_{N_{\xi }}=\frac{1}{2}\int \left( \varepsilon _{2}+\varepsilon
_{1}\left( \ Curl \overrightarrow{N}.\overrightarrow{B}\right)
^{2}+\varepsilon _{3}\left( \ Div \overrightarrow{B}\right) ^{2}\right)
d\xi .  \tag{5.5}
\end{equation}

Similarly, since $\nabla _{B}B=-\Psi _{B}\overrightarrow{T}+-\varepsilon
_{2}\left( \ Div \overrightarrow{B}\right) \overrightarrow{N}$ the
energy of the vector field $\frac{\partial B}{\partial \xi }$ is expressed
as 
\begin{equation}
energy_{B_{\xi }}=\frac{1}{2}\int \left( \varepsilon _{3}+\varepsilon
_{1}\left( \ Curl \overrightarrow{B}.\overrightarrow{B}\right)
^{2}+\varepsilon _{2}\left( \ Div \overrightarrow{B}\right) ^{2}\right)
d\xi .  \tag{5.6}
\end{equation}

\subsection{The energy of the tangent vector of $\protect\eta -$lines on a
moving particle in $M_{q}^{3}\left( c\right) $}

In the subsection, the bending energy formulas of the unit tangent vector of 
$\eta -$lines an elastic curve given by extended Serret-Frenet relations
along the curve $\gamma (s,\xi ,\eta )$ are expressed in $M_{q}^{3}\left(
c\right) .$ For the curve $\gamma (s,\xi ,\eta )$ with parameter $\eta $,
which $\eta $ is the distance along the $\eta -$lines of the curve in $\eta -
$direction and tangent vector of $\eta -$lines is described by $\frac{%
\partial T}{\partial \eta },$ from Sasaki metric the energy on the particle
in vector field $\frac{\partial T}{\partial \eta }$ is written as 
\begin{equation*}
energy_{T_{\eta }}=\frac{1}{2}\int \rho _{\eta }(dT(T),dT(T))d\eta ,
\end{equation*}%
from (2.3), (2.4), (2.5) and we get%
\begin{equation*}
\rho _{\eta }(dT(T),dT(T))=\rho _{\eta }\left( T,T\right) +\rho _{\eta
}\left( \nabla _{T}T,\nabla _{T}T\right) 
\end{equation*}%
also from extended Serret-Frenet relations with respect to parameter $\eta $
or since $\nabla _{T}T=\varepsilon _{1}\left( -\varepsilon _{3}\tau -\Psi
_{N}\right) \overrightarrow{N}-\varepsilon _{1}\varepsilon _{3}(\varepsilon
_{2}\tau +\ Curl \overrightarrow{B}\cdot \overrightarrow{N})%
\overrightarrow{B}$, we get 
\begin{equation}
energy_{T_{\eta }}=\frac{1}{2}\int (\varepsilon _{1}+\varepsilon
_{2}(\varepsilon _{3}\tau +\ Curl \overrightarrow{N}\cdot 
\overrightarrow{N})^{2}+\varepsilon _{3}(\varepsilon _{2}\tau +\ Curl%
\overrightarrow{B}\cdot \overrightarrow{N})^{2})d\eta .  \tag{5.7}
\end{equation}

Similarly, since $\nabla _{N}N=\left( \varepsilon _{3}\tau +\Psi _{N}\right) 
\overrightarrow{T}+\varepsilon _{2}\varepsilon _{3}\ Curl%
\overrightarrow{B}\cdot \overrightarrow{T}\overrightarrow{B}$ $,$ the energy
of the vector field $\frac{\partial N}{\partial \eta }$ is written as 
\begin{equation}
energy_{N_{\eta }}=\int \left( \varepsilon _{2}+\varepsilon _{1}(\varepsilon
_{3}\tau +\ Curl\overrightarrow{N}\cdot \overrightarrow{N}%
)^{2}+\varepsilon _{2}(\ Curl \overrightarrow{B}\cdot \overrightarrow{T}%
)^{2}\right) d\eta .  \tag{5.8}
\end{equation}%
and since $\nabla _{B}B=(\varepsilon _{2}\tau +\ Curl \overrightarrow{B}%
\cdot \overrightarrow{N})\overrightarrow{T}+(-\ Curl \overrightarrow{B}%
\cdot \overrightarrow{T})\overrightarrow{N},$ the energy of the vector field 
$\frac{\partial B}{\partial \eta }$ is also obtained as 
\begin{equation}
energy_{B_{\eta }}=\frac{1}{2}\int \left( \varepsilon _{3}+\varepsilon
_{1}(\varepsilon _{2}\tau +\ Curl \overrightarrow{B}\cdot 
\overrightarrow{N})^{2}+\varepsilon _{2}(\ Curl \overrightarrow{B}\cdot 
\overrightarrow{T})^{2}\right) d\eta .  \tag{5.9}
\end{equation}

\end{document}